\documentclass[11pt]{amsart}

\usepackage{amscd,array,graphicx,color,epsfig}

\newtheorem{theorem}{Theorem}[section]
\newtheorem{proposition}[theorem]{Proposition}
\newtheorem{lemma}[theorem]{Lemma}
\newtheorem{corollary}[theorem]{Corollary}

\newtheorem*{parameterizationTheorem}{Parameterization Theorem}
\newtheorem*{MinimumBridgeNumberTheorem}{Minimum Bridge Number Theorem~\ref{thm:minbridgenum}}
\newtheorem{minimumBridgeNumberTheorem}[theorem]{Minimum Bridge Number Theorem}
\newtheorem*{ubtheorem}{Theorem~\ref{thm:ub}}
\newtheorem*{maxbridge}{Theorem~\ref{thm:max}}

\theoremstyle{definition}
\newtheorem{definition}[theorem]{Definition}

\newtheorem{remark}[theorem]{Remark}

\numberwithin{equation}{section}

\newcommand{\bc}{\operatorname{bc}}
\newcommand{\br}{\operatorname{br}}

\newcommand{\depth}{\operatorname{depth}}
\newcommand{\dist}{\operatorname{dist}}
\newcommand{\lk}{\operatorname{lk}}
\newcommand{\Nbd}{\operatorname{Nbd}}

\newcommand{\rbc}{\operatorname{rbc}}

\renewcommand{\P}{\operatorname{{\mathcal P}}}
\newcommand{\Q}{\operatorname{{\mathbb Q}}}
\newcommand{\R}{\operatorname{{\mathbb R}}}
\newcommand{\Z}{\operatorname{{\mathbb Z}}}

\newcommand{\D}{\operatorname{{\mathcal D}}}
\newcommand{\G}{\operatorname{{\mathcal G}}}
\newcommand{\K}{\operatorname{{\mathcal K}}}
\newcommand{\T}{\operatorname{{\mathcal T}}}

\begin{document}

\title[The depth of a knot tunnel]
{The depth of a knot tunnel}

\author{Sangbum Cho}
\address{Department of Mathematics\\
University of Oklahoma\\
Norman, Oklahoma 73019\\
USA} 
\email{scho@ou.edu}

\author{Darryl McCullough}
\address{Department of Mathematics\\
University of Oklahoma\\
Norman, Oklahoma 73019\\
USA} 
\email{dmccullough@math.ou.edu}
\urladdr{www.math.ou.edu/$_{\widetilde{\phantom{n}}}$dmccullough/}


\date{\today}


\begin{abstract}
The theory of tunnel number $1$ knots detailed in \cite{CM} provides a
non-negative integer invariant $\depth(\tau)$ for a knot tunnel $\tau$. We
give various results related to the depth invariant. Noting that it equals
the minimal number of Goda-Scharlemann-Thompson ``tunnel moves'' \cite{GST}
needed to construct the tunnel, we calculate the number of distinct minimal
sequences of tunnel moves that can produce a given tunnel. Next, we give a
recursion that tells the minimum bridge number of a knot having a tunnel of
depth $d$. The growth of this value is proportional to $(1+\sqrt{2})^d$,
which improves known estimates of the rate of growth of bridge number as a
function of the Hempel distance of the associated Heegaard splitting. We
also give various upper bounds for bridge number in terms of the cabling
constructions needed to produce a tunnel of a knot, showing in particular
that the maximum bridge number of a knot produced by $n$ cabling
constructions is the $(n+2)^{nd}$ Fibonacci number.  Finally, we explicitly
compute the slope parameters for the regular (or ``short'') tunnels of
torus knots, and find a sequence of them for which the bridge numbers of
the associated knots achieve the growth rate~$(1+\sqrt{2})^d$.\par
\end{abstract}

\maketitle

\section*{Introduction}
\label{sec:intro}

This work concerns a new invariant of knot tunnels, called the
\textit{depth}. It is based on the theory of knot tunnels developed in our
earlier work \cite{CM}, which provides a simplicial complex $\D(H)/\Gamma$
whose vertices correspond to the (equivalence classes of) tunnels of all
tunnel number~$1$ knots. The depth invariant of a tunnel is defined to be
the simplicial distance in the $1$-skeleton of $\D(H)/\Gamma$ from the
vertex corresponding to the tunnel to the vertex corresponding to the
unique tunnel of the trivial knot. In particular, the trivial tunnel is the
only tunnel of depth~$0$. A tunnel has depth~$1$ exactly when it is a
$(1,1)$-tunnel of a $(1,1)$-knot.

We denote the depth of a tunnel $\tau$ by $\depth(\tau)$. It is somewhat
similar to the \textit{$($Hempel\/$)$ distance} $\dist(\tau)$ (see
J. Johnson \cite{JohnsonBridgeNumber} and Y. Minsky, Y. Moriah, and
S. Schleimer \cite{MMS}), but is very easy to calculate in terms of the
parameter description of tunnels given in~\cite{CM}.
The two invariants are related by the inequality
\[\dist(\tau)-1 \leq \depth(\tau)\ ,\]
but the depth can be much larger than the distance. Indeed, we will see
that the regular tunnels of torus knots have distance~$2$, but their depths
can be arbitrarily large.

The depth invariant has a geometric interpretation in terms of a
construction that first appeared in a paper of H. Goda, M. Scharlemann, and
A. Thompson \cite{GST}. Their construction, which we call a giant step, takes
a tunnel and produces a new tunnel (usually of a different knot). They
proved that every tunnel could be produced starting from the tunnel of the
trivial knot and applying a sequence of giant steps, and we will see from the
definitions that $\depth(\tau)$ is the minimum length of a such a sequence.
Unlike the construction of a knot tunnel using cabling operations,
developed in~\cite{CM} and reviewed in section~\ref{sec:cabling} below, the
choice of giant steps is usually not unique, even when one restricts to
minimal sequences. Using the simplicial structure of $\D(H)/\Gamma$, we
will give an algorithm to calculate the number of distinct minimal
sequences of giant steps that produce a given tunnel. In particular, this
provides arbitrarily complicated examples of tunnels for which the minimal
giant steps sequence \textit{is} unique, while showing that such tunnels are
sparse among the set of all tunnels. The algorithm is quite effective. We
have implemented it computationally~\cite{slopes} to find the number of
distinct minimal sequences producing a tunnel, given its parameter
description from~\cite{CM}.

We next turn to an examination of the bridge number of a tunnel number~$1$
knot. The first main result, theorem~\ref{thm:bridge_numbers}, gives a
lower bound for the bridge number of a specific tunnel number~$1$ knot in
terms of the parameters of a tunnel of the knot. The algorithm to compute
that bound is easy, and we have implemented it
computationally~\cite{slopes}. The proof of
theorem~\ref{thm:bridge_numbers} is quite easy for us since the necessary
ideas and hard geometric work were already developed by Goda, Scharlemann,
and Thompson \cite{GST} and Scharlemann and Thompson
\cite{Scharlemann-Thompson}.  Theorem~\ref{thm:bridge_numbers}, together
with a geometric construction involving the cabling construction of tunnels
from~\cite{CM}, gives a general and sharp lower bound for the bridge number
in terms of the depth of a tunnel:
\begin{MinimumBridgeNumberTheorem} 
For $d\geq 1$, the minimum bridge number of a knot having a tunnel of
depth~$d$ is $a_d$, where $a_1=2$, $a_2=4$, and $a_d=2a_{d-1}+a_{d-2}$ for
$d\geq 3$.
\end{MinimumBridgeNumberTheorem}
\noindent As a matrix, the recursion in theorem~\ref{thm:minbridgenum} is
\[\begin{pmatrix} a_{d+1}\\
a_d\end{pmatrix} =
\begin{pmatrix} 2 & 1 \\
1 & 0
\end{pmatrix}
\begin{pmatrix} a_d\\
a_{d-1}
\end{pmatrix}\ .
\]
The eigenvalues of this matrix are $1\pm \sqrt{2}$, showing that the
asymptotic growth rate of the bridge numbers of any sequence of tunnel
number~$1$ knots as a function of depth is at least a constant multiple of
$(1+\sqrt{2})^d$. This improves Lemma~2 of \cite{JohnsonBridgeNumber},
which is that bridge number grows linearly with distance. Since each of the
giant steps of a minimal sequence increases the depth by~$1$,
corollary~\ref{coro:bridge_numbers} also improves Proposition~1.11 of
\cite{GST}, which implies that bridge number grows asymptotically at least
as fast as~$2^d$.

Of course, the Minimum Bridge Number Theorem also shows that the bound for
growth rate of $(1+\sqrt{2})^d$ is best possible, indeed its proof tells
exactly how to construct a tunnel of depth~$d$ having bridge
number~$a_d$. We will show that this growth rate is also achieved by a
sequence of torus knot tunnels (each obtained by applying a cabling
operation to the previous one) given in section~\ref{sec:torus_knots}. In
fact, the bridge numbers of the knots of that sequence are given by the
recursion in the Minimum Bridge Number Theorem, except that one starts with
$a_1=2$ and $a_2=5$. The terms $a_d$ are then the minimal bridge numbers of
any torus knot having a tunnel of the corresponding depth.

Besides the specific examples of torus knot tunnels we have already
mentioned, we give a general algorithm to compute the slope parameters for
the regular tunnel (sometimes called the ``short'' tunnel) of a $(p,q)$
torus knot. The algorithm uses the continued fraction expansion of
$p/q$. It is very effective and has been implemented
computationally~\cite{slopes}.

A more general version of the geometric construction used in proving the
Minimum Bridge Number Theorem allows us to show that in general, cabling
operations can be carried out rather efficiently with respect to bridge
number. This leads to various upper bounds for the bridge number of a knot
in terms of its tunnels. In particular,
\begin{ubtheorem} 
Let $(F_1, F_2, \ldots)=(1,1,2,3,\ldots)$ be the Fibonacci sequence.
Suppose that $\tau$ is a regular tunnel produced by $n$ cabling operations,
of which the first $m$ produce semisimple tunnels. Then
$\br(K_{\tau_n})\leq mF_{n-m+2}+F_{n-m+1}$.
\end{ubtheorem}
\noindent For fixed $n$, the largest value for the upper bound in
theorem~\ref{thm:ub} occurs when $m=2$, and we show that it is sharp
for this case:
\begin{maxbridge} The maximum bridge number of any tunnel number
$1$-knot having a tunnel produced by $n$ cabling operations is $F_{n+2}$.
\end{maxbridge}

Here is an outline of the sections of the paper. The first three sections
constitute a concise review of material from~\cite{CM} that we will need
for the present applications. Section~\ref{sec:ddd} introduces the distance
and depth invariants, and gives a few results that follow quickly
from~\cite{CM} and work of other authors. The main applications of the
paper may be then read independently. The giant steps discussed above are
introduced in section~\ref{sec:GST_moves}, and the analysis of minimal
sequences of giant steps is carried out in
section~\ref{sec:GST_move_seqs}. Lower bounds for bridge number, in
particular the Minimum Bridge Number Theorem, are given in
section~\ref{sec:bridge_number_growth}, while
section~\ref{sec:bridgenum_conj} has the results on upper bounds.  The
torus knot examples are worked out in section~\ref{sec:torus_knots}. The
final section of the paper reviews how the general theory can be adapted to
include tunnel number~$1$ links, and indicates how the applications in the
present paper extend to that case.

\section{The disk complex of an irreducible $3$-manifold}
\label{sec:disk_complex}

\begin{figure}
\begin{center}
\includegraphics[width=28 ex]{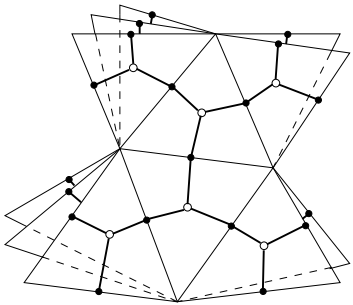}
\caption{A portion of the nonseparating disk complex $\D(H)$ and the tree
$\widetilde{\T}$. Countably many $2$-simplices meet along each edge.}
\label{fig:subdivision}
\end{center}
\end{figure}

Let $H$ be a genus~$2$ orientable handlebody, regarded as the standard
unknotted handlebody in $S^3$. For us, a \textit{disk in H} means a
properly imbedded disk in $H$, \textit{which is assumed to be nonseparating
unless otherwise stated.} The \textit{disk complex} $\D(H)$ is a
$2$-dimensional, contractible simplicial complex, whose vertices are the
(proper) isotopy classes of essential properly imbedded disks in $H$, such
that a collection of $k+1$ vertices spans a $k$-simplex if and only if they
admit a set of pairwise-disjoint representatives. Each $1$-simplex of
$\D(H)$ is a face of countably many $2$-simplices. As suggested by
figure~\ref{fig:subdivision}, $\D(H)$ grows outward from any of its
$2$-simplices in a treelike way. In fact, it deformation retracts to the
tree $\widetilde{\T}$ seen in figure~\ref{fig:subdivision}.

Each tunnel of a tunnel number 1 knot determines a collection of disks in
$H$ as follows. The tunnel is a $1$-handle attached to a regular
neighborhood of the knot to form an unknotted genus-$2$ handlebody. An
isotopy carrying this handlebody to $H$ carries a cocore $2$-disk of that
$1$-handle to a nonseparating disk in $H$, and carries the tunnel
number~$1$ knot to a core circle of the solid torus obtained by cutting $H$
along the image disk in $H$. The indeterminacy of this isotopy is the
group of isotopy classes of orientation-preserving homeomorphisms of $S^3$
that preserve $H$. This group is called the \textit{Goeritz group $\G$.}
Work of M. Scharlemann \cite{ScharlemannTree} and E. Akbas \cite{Akbas}
proves that $\G$ is finitely presented, and even provides a simple
presentation of it.

\begin{figure}\begin{center}
\includegraphics[width=28 ex]{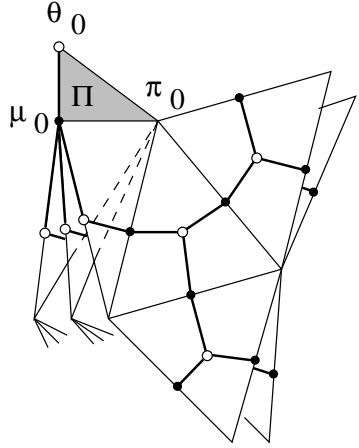}
\caption{A portion of $\D(H)/\G$ and $\T$ near the primitive orbits.}
\label{fig:Delta}
\end{center}
\end{figure}

Since two disks in $H$ determine equivalent tunnels exactly when they
differ by an isotopy moving $H$ through $S^3$, \textit{the collection of
all tunnels of all tunnel number~$1$ knots corresponds to the set of orbits
of vertices of $\D(H)$ under $\G$.} So it is natural to examine the
quotient complex $\D(H)/\G$, which is illustrated in figure~\ref{fig:Delta}.
Through work of the first author~\cite{Cho}, this action is
well-understood. A \textit{primitive} disk in $H$ is a disk $D$ such that
there is a disk $E$ in $\overline{S^3-H}$ for which $\partial D$ and
$\partial E$ intersect transversely in one point in $\partial H$.  The
primitive disks (regarded as vertices) span a contractible subcomplex
$\P(H)$ of $\D(H)$, called the \textit{primitive subcomplex}. The action of
$\G$ on $\P(H)$ is as transitive as possible, indeed the quotient
$\P(H)/\G$ is a single $2$-simplex $\Pi$ which is the image of any
$2$-simplex of the first barycentric subdivision of $\P(H)$. Its vertices
are $\pi_0$, the orbit of all primitive disks, $\mu_0$, the orbit of all
pairs of disjoint primitive disks, and $\theta_0$, the orbit of all triples
of disjoint primitive disks.

On the remainder of $\D(H)$, the stabilizers of the action are as small as
possible. A $2$-simplex which has two primitive vertices and one
nonprimitive is identified with some other such simplices, then folded in
half and attached to $\Pi$ along the edge $\langle \mu_0,\pi_0\rangle$.
The nonprimitive vertices of such $2$-simplices are exactly the disks in
$\D(H)$ that are disjoint from some primitive pair, and these are called
$\textit{simple}$ disks. As tunnels, they are the upper and lower tunnels
of $2$-bridge knots. The remaining $2$-simplices of $\D(H)$ receive no
self-identifications, and descend to portions of $\D(H)/\G$ that are
treelike and are attached to one of the edges $\langle \pi_0,\tau_0\rangle$
where $\tau_0$ is simple.

The tree $\widetilde{\T}$ shown in figure~\ref{fig:subdivision} is
constructed as follows.  Let $\D'(H)$ be the first barycentric subdivision
of $\D(H)$. Denote by $\widetilde{\T}$ the subcomplex of $\D'(H)$ obtained
by removing the open stars of the vertices of $\D(H)$. It is a bipartite
graph, with ``white'' vertices of valence $3$ represented by triples and
``black'' vertices of (countably) infinite valence represented by
pairs. The valences reflect the fact that moving along an edge from a
triple to a pair corresponds to removing one of its three disks, while
moving from a pair to a triple corresponds to adding one of infinitely many
possible third disks to a pair.  The possible disjoint third disks that can
be added are called the \textit{slope disks} for the pair.

The image $\widetilde{\T}/\G$ of $\widetilde{\T}$ in $\D(H)/\G$ is a tree
$\T$.  The vertices of $\D(H)/\G$ that are images of vertices of $\D(H)$
are not in $\T$, but their links in $\D'(H)/\G$ are subcomplexes of
$\T$. These links are infinite trees. For each such vertex $\tau$ of
$\D(H)/\G$, i.~e.~each tunnel, there is a unique shortest path in $\T$ from
$\theta_0$ to \textit{the vertex in the link of $\tau$ that is closest to
$\theta_0$.}  This path is called the \textit{principal path} of $\tau$,
and this closest vertex is a triple, called the \textit{principal vertex}
of $\tau$. The two disks in the principal vertex, other than $\tau$, are
called the \textit{principal meridian pair} of $\tau$. They are exactly the
disks called $\mu^+$ and $\mu^-$ that play a key role
in~\cite{Scharlemann-Thompson}. Figure~\ref{fig:corridor} below shows the
principal path of a certain tunnel.

\section{Slope disks}
\label{sec:slope_disks}

In \cite{CM}, it is explained how moving along the principal path of a
tunnel $\tau$ encodes a sequence of cabling constructions starting from the
tunnel of the trivial knot and producing a sequence of tunnels ending with
$\tau$. We will review the cabling construction in
section~\ref{sec:cabling} below. Each cabling is determined by a rational
parameter, called its \textit{slope}. In section~\ref{sec:torus_knots}
below, we will compute these slopes for many tunnels of torus knots, so it
is necessary to recall the precise definition of slope. For the other
applications, the precise details of the definition are not needed, so a
more superficial reading of this section might suffice.

\begin{figure}
\begin{center}
\includegraphics[width=65 ex]{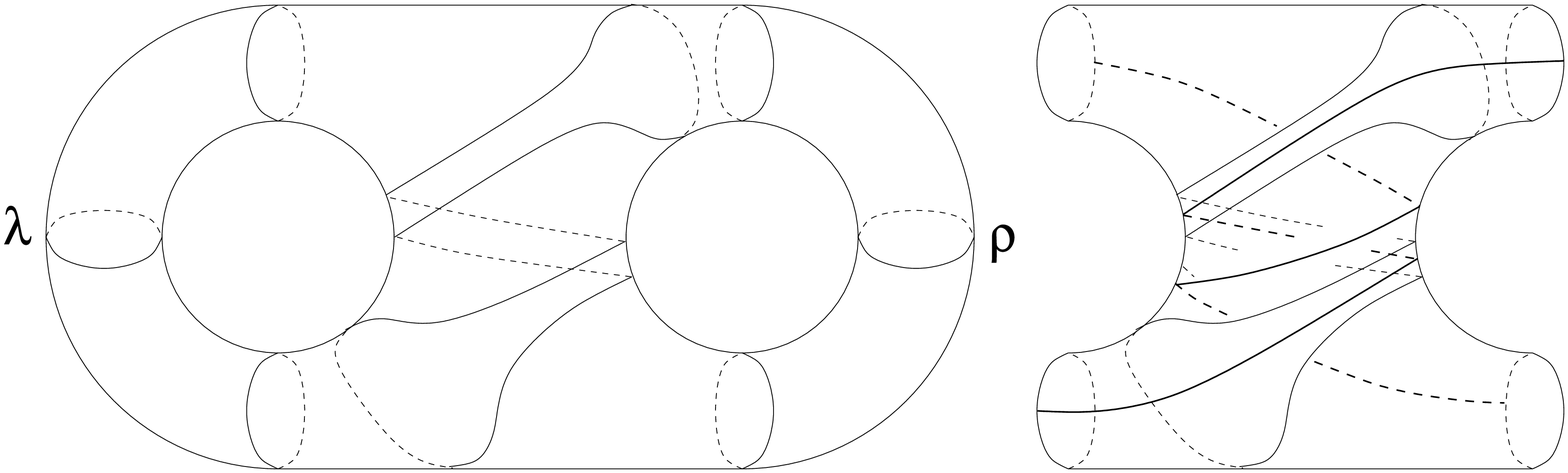}
\caption{A slope disk and one of its cabling arcs.}
\label{fig:slopes}
\end{center}
\end{figure}

Fix a pair of disks $\lambda$ and $\rho$ (for ``left'' and ``right'') in
$H$, as shown abstractly in figure~\ref{fig:slopes}. Of course, in the true
picture in $H$ in $S^3$, these can look a great deal more complicated than
the primitive pair shown in figure~\ref{fig:slopes}. Let $B$ be $H$ cut
along $\lambda\cup \rho$. The frontier of $B$ in $H$ consists of four disks
which appear vertical in figure~\ref{fig:slopes}. Denote this frontier by
$F$, and let $\Sigma$ be $B\cap \partial H$, a sphere with four holes.
A \textit{slope disk for $\{\lambda,\rho\}$} is an essential disk, possibly
separating, which is contained in $B$ and not isotopic to any component
of~$F$. The boundary of a slope disk always separates $\Sigma$ into two
pairs of pants, conversely any loop in $\Sigma$ that is not homotopic into
$\partial \Sigma$ is the boundary of a unique slope disk. If two slope
disks are isotopic in $H$, then they are isotopic in~$B$.

An arc in $\Sigma$ whose endpoints lie in two different boundary circles of
$\Sigma$ is called a \textit{cabling arc.}  Figure~\ref{fig:slopes} shows a
pair of cabling arcs disjoint from a slope disk. A slope disk is disjoint
from a unique pair of cabling arcs, and each cabling arc determines a
unique slope disk.

\begin{figure}
\begin{center}
\includegraphics[width=45 ex]{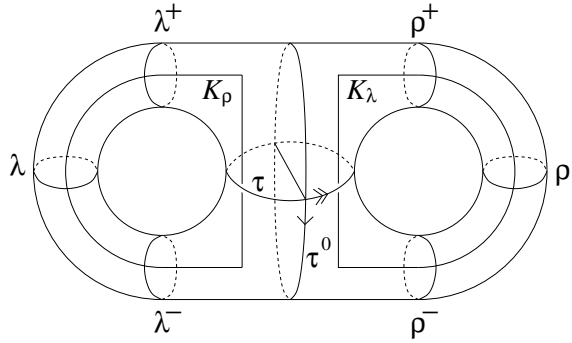}
\caption{The slope-zero perpendicular disk $\tau^0$. It is chosen so that
$K_\lambda$ and $K_\rho$ have linking number~$0$.}
\label{fig:slope_coords}
\end{center}
\end{figure}

Each choice of nonseparating slope disk for a pair $\mu=\{\lambda,\rho\}$
determines a correspondence between~$\Q\cup\{\infty\}$ and the set of all
slope disks of $\mu$, as follows. Fixing a nonseparating slope disk $\tau$
for $\mu$, write $(\mu;\tau)$ for the ordered pair consisting of $\mu$ and
$\tau$.
\begin{definition} A \textit{perpendicular disk} for $(\mu;\tau)$
is a disk $\tau^\perp$, with the following properties:
\begin{enumerate}
\item $\tau^\perp$ is a slope disk for $\mu$.
\item $\tau$ and $\tau^\perp$ intersect transversely in one arc.
\item $\tau^\perp$ separates $H$.
\end{enumerate}
\end{definition}
There are infinitely many choices for $\tau^\perp$, but because $H\subset
S^3$ there is a natural way to choose a particular one, which we call
$\tau^0$. It is illustrated in figure~\ref{fig:slope_coords}. To construct
it, start with any perpendicular disk and change it by Dehn twists of $H$
about $\tau$ until the core circles of the complementary solid tori have
linking number~$0$ in~$S^3$.

For calculations, it is convenient to draw the picture as in
figure~\ref{fig:slope_coords}, and orient the boundaries of $\tau$ and
$\tau^0$ so that the orientation of $\tau^0$ (the ``$x$-axis''), followed
by the orientation of $\tau$ (the ``$y$-axis''), followed by the outward
normal of $H$, is a right-hand orientation of $S^3$. At the other
intersection point, these give the left-hand orientation, but the
coordinates are unaffected by changing the choices of which of
$\{\lambda,\rho\}$ is $\lambda$ and which is $\rho$, or changing which of
the disks $\lambda^+$, $\lambda^-$, $\rho^+$, and $\rho^-$ are ``$+$'' and
which are ``$-$'', provided that the ``$+$'' disks both lie on the same
side of $\lambda\cup\rho\cup\tau$ in figure~\ref{fig:slope_coords}.

\begin{figure}
\begin{center}
\includegraphics[width=\textwidth]{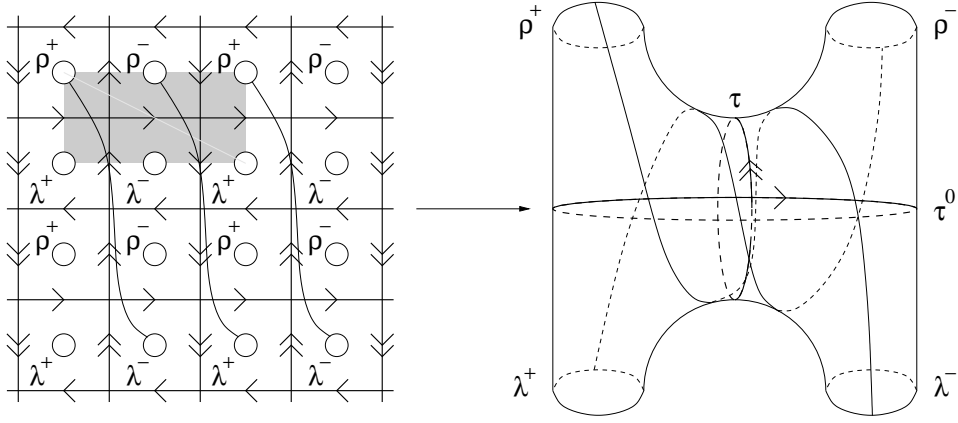}
\caption{The covering space $\widetilde{\Sigma}\to\Sigma$, and some lifts
of a $[1,-3]$-cabling arc. The shaded region is a fundamental domain.}
\label{fig:covering}
\end{center}
\end{figure}

Let $\widetilde{\Sigma}$ be the covering space of $\Sigma$ such that:
\begin{enumerate}
\item $\widetilde{\Sigma}$ is the plane with an open disk of radius $1/8$
removed from each point with half-integer coordinates.
\item The components of the preimage of $\tau$ are the vertical lines
with integer $x$-coordinate.
\item The components of the preimage of $\tau^0$ are the horizontal lines
with integer $y$-coordinate.
\end{enumerate}
\noindent Figure~\ref{fig:covering} shows a picture of $\widetilde{\Sigma}$
and a fundamental domain for the action of its group of covering
transformations, which is the orientation-preserving subgroup of the group
generated by reflections in the half-integer lattice lines (that pass
through the centers of the missing disks). Each circle of
$\partial\widetilde{\Sigma}$ double covers a circle of~$\partial \Sigma$.

If we lift any cabling arc in $\Sigma$ to $\widetilde{\Sigma}$, the lift
runs from a boundary circle of $\widetilde{\Sigma}$ to one of its
translates by a vector $(p,q)$ of signed integers, defined up to
multiplication by the scalar $-1$. Thus each cabling arc receives a
\textit{slope pair} $[p,q]=\{(p,q),(-p,-q)\}$, and is called a
\textit{$[p,q]$-cabling arc.} The corresponding slope disk receives the
slope pair $[p,q]$ as well.

An important observation is that a $[p,q]$-slope disk is nonseparating in
$H$ if and only if $q$ is odd. Both happen exactly when the corresponding
cabling arc has one endpoint in $\lambda^+$ or $\lambda^-$ and the other in
$\rho^+$ or~$\rho^-$.

\begin{definition} The \textit{$(\mu;\tau)$-slope} of a $[p,q]$-slope 
disk or cabling arc is~$q/p\in \Q\cup\{\infty\}$.
\end{definition}
\noindent The $(\mu;\tau)$-slope of $\tau^0$ is $0$, and the
$(\mu;\tau)$-slope of $\tau$ is~$\infty$.

Slope disks for a \textit{primitive} pair are handled in a special way.
Rather than using a particular choice of $\tau$ from the context, one
chooses $\tau$ to be some third primitive disk. Altering this choice can
change $[p,q]$ to any $[p+nq,q]$, but the quotient $p/q$ is well-defined as
an element of $\Q/\Z\cup\{\infty\}$. This element $[p/q]$ is called the
\textit{simple slope} of the slope disk (it is $[0]$ exactly when the slope
disk is itself primitive). Two simple disks have the same simple slope
exactly when they are equivalent by an element of the Goeritz group.

\section{Parameterization and cabling operations}
\label{sec:cabling}

In this section, we discuss the Parameterization Theorem
from~\cite{CM}. First, we review the cabling construction.

In a sentence, the cabling construction is to ``Think of the union of $K$
and the tunnel arc as a $\theta$-curve, and rationally tangle the ends of
the tunnel arc and one of the arcs of $K$ in a neighborhood of the other
arc of $K$.''  We sometimes call this ``swap and tangle,'' since one of the
arcs in the knot is exchanged for the tunnel arc, then the ends of other
arc of the knot and the tunnel arc are connected by a rational tangle.

\begin{figure}
\begin{center}
\includegraphics[width=\textwidth]{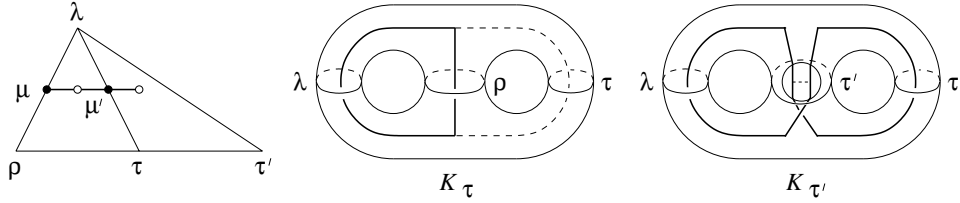}
\caption{Schematic for the general cabling construction. In the middle ball
in the right-hand picture of $H$, the two vertical arcs form some rational
tangle, disjoint from the disk~$\tau'$.}
\label{fig:schematic}
\end{center}
\end{figure}
Figure~\ref{fig:schematic} illustrates the cabling construction
schematically.  Begin with a pair $\mu=\{\lambda,\rho\}$ and a triple
$\mu\cup \{\tau\}$. In a $\theta$-curve corresponding to $\mu\cup
\{\tau\}$, the union of the arcs dual to $\mu$ is $K_\tau$, and the arc
dual to $\tau$ is a tunnel arc for $K_\tau$. Moving through $\T$ starting
at the edge $\langle \mu,\mu\cup\{\tau\}\rangle$ determines a sequence of
steps in which one of the two disks of a pair $\{\lambda,\rho\}$ is
replaced by a tunnel disk $\tau$, and a slope disk $\tau'$ of the new pair
$\mu'$ (with $\tau'$ nonseparating in $H$) is chosen as the new tunnel
disk, ending up at the edge $\langle \mu',\mu'\cup\{\tau'\}\rangle$. This
is a cabling operation producing $\tau'$ from $\tau$. It is required that
$\tau'\neq \tau$, that is, cablings do not allow one to ``backtrack'' in
$\T$.

As illustrated in figure~\ref{fig:schematic}, the way that the path
determines the particular cabling operation is:
\begin{enumerate}
\item The selection of $\lambda$ or $\rho$ corresponds to which edge one
chooses to move out of the white vertex $\{\lambda,\rho,\tau\}$.
\item The selection of the new slope disk $\tau'$ corresponds to which edge
one chooses to continue  out of the black vertex~$\mu'$.
\end{enumerate}

\begin{figure}
\begin{center}
\includegraphics[width=\textwidth]{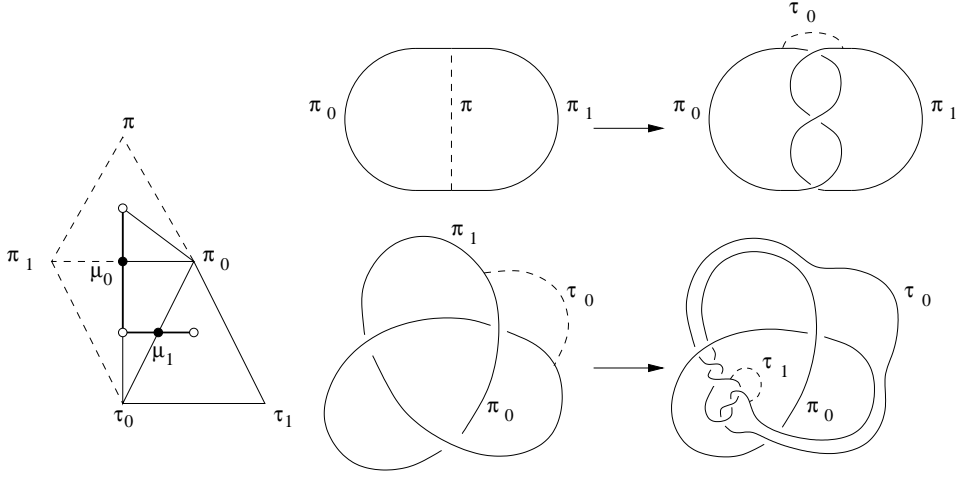}
\caption{Examples of the cabling construction.}
\label{fig:cabling}
\end{center}
\end{figure}

Figure~\ref{fig:cabling} shows the effects of a specific sequence of two
cabling constructions, starting with the trivial knot and obtaining the
trefoil, then a cabling construction starting with the tunnel of the
trefoil.

When $\mu$ is the primitive pair $\mu_0$, and $\tau_0$ is a simple disk for
$\mu_0$, the pair $(\mu_0;\tau_0)$ determines a cabling construction
starting with the tunnel of the trivial knot and producing $\tau_0$, which
is an upper or lower tunnel of a $2$-bridge knot. This is a \textit{simple
cabling of slope $m_0$,} where $m_0$ is the simple slope of~$\tau_0$.

The principal path of a tunnel $\tau$ determines a sequence of cablings,
which produce $\tau$ from the trivial tunnel $\pi_0$. The first is a simple
cabling of some simple slope $m_0\in \Q/\Z$. For $i\geq 1$, the cabling
producing $\tau_i$ from $\tau_{i-1}$ obtains a rational slope $m_i$ as
follows.  Let $\sigma$ be the unique disk of $\mu_{i-1}-\mu_i$ (the
``trailing'' disk), which in the schematic picture of
figure~\ref{fig:schematic} happens to be $\rho$. The slope of this cabling
is defined to be the rational number $m_i$ which is the
$(\mu_i;\sigma)$-slope of $\tau_i$. Note that an illegal ``backtrack''
cabling would have slope~$\infty$, since the $(\mu;\sigma)$-slope pair of
$\sigma$ is~$[0,1]$.

The Unique Cabling Sequence Theorem, theorem~13.2 of~\cite{CM}, states that
there is a unique sequence of cabling constructions starting with the
tunnel of the trivial knot and ending with $\tau$. It is an immediate
consequence of the fact that $\T$ is a tree. Viewed in terms of slope, this
becomes the following result, theorem~12.3 of~\cite{CM}:
\begin{parameterizationTheorem}
Let $\tau$ be a knot tunnel with principal path $\theta_0$, $\mu_0$,
$\mu_0\cup \{\tau_0\}$, $\mu_1,\ldots\,$, $\mu_n$, $\mu_n\cup \{\tau_n\}$.
Fix a lift of the principal path to $\D(H)$, so that each $\mu_i$
corresponds to an actual pair of disks in~$H$.
\begin{enumerate}
\item If $\tau$ is primitive, put $m_0=[0]\in\Q/\Z$. Otherwise,
let $m_0=[p_0/q_0]\in\Q/\Z$ be the simple slope of $\tau_0$.
\item If $n\geq 1$, then for $1\leq i\leq n$ let $\sigma_i$ be the unique disk
in $\mu_{i-1}-\mu_i$ and let $m_i=q_i/p_i\in\Q$ be the
$(\mu_i;\sigma_i)$-slope of $\tau_i$.
\item If $n\geq 2$, then for $2\leq i\leq n$ define $s_i=0$ or $s_i=1$
according to whether or not the unique disk of $\mu_i\cap\mu_{i-1}$ equals the
unique disk of $\mu_{i-1}\cap\mu_{i-2}$.
\end{enumerate}
Then, sending $\tau$ to the pair $((m_0,\ldots,m_n),(s_2,\ldots,s_n))$ is a
bijection from the set of all tunnels of all tunnel number~$1$ knots to the
set of all elements $(([p_0/q_0],q_1/p_1,\ldots,q_n/p_n),(s_2,\ldots,s_n))$
in
\[\big(\Q/\Z \big)\,\cup\, 
\big(\Q/\Z\,\times\, \Q\big) \,\cup\, \big(\cup_{n\geq 2}
\;\Q/\Z\,\times\, \Q^n \,\times\,\, C_2^{n-1}\big)\]
with all $q_i$ odd.
\end{parameterizationTheorem}

\begin{remark} Up to details of definition, the final slope $m_n$ in the
Parameterization Theorem is the Scharlemann-Thompson
invariant~\cite{Scharlemann-Thompson}.
\end{remark}

A tunnel produced from the tunnel of the trivial knot by a single cabling
construction is called a \textit{simple} tunnel. As already noted, these
are the ``upper and lower'' tunnels of $2$-bridge knots. According to the
Parameterization Theorem, these are determined by a single $\Q/\Z$-valued
parameter $m_0$, and this is of course a version of the standard rational
parameter that classifies the $2$-bridge knot.

A tunnel is called \textit{semisimple} if it is disjoint from a primitive
disk, but not from any primitive pair. The simple and semisimple tunnels
are exactly the $(1,1)$-tunnels, that is, the tunnels that can be put into
$1$-bridge position with respect to a Heegard torus of $S^3$. The non-simple
tunnels of $2$-bridge knots are semisimple, and in~\cite{CM}, their
parameter sequences are calculated. Since they are semisimple, their
parameters $s_2,\ldots\,$, $s_n$ in the Parameterization Theorem are all
$0$. Their slope parameters are determined by a somewhat complicated, but
easily programmable algorithm using the continued fraction expansion of the
classifying parameter.

Finally, a tunnel is called \textit{regular} if it is neither primitive,
simple, or semisimple.

\section{Distance and depth}
\label{sec:ddd}

In this section we formally introduce the distance and depth invariants of
a tunnel~$\tau$. The \textit{$($Hempel\/$)$ distance} $\dist(\tau)$ is the
shortest distance in the curve complex of $\partial H$ from $\partial \tau$
to a loop that bounds a disk in $\overline{S^3-H}$ (see J. Johnson
\cite{JohnsonBridgeNumber} and Y. Minsky, Y. Moriah, and S. Schleimer
\cite{MMS}). It is well-defined since the action of the Goeritz group on
$\partial H$ preserves the set of loops that bound disks in $H$ and the set
that bound in~$\overline{S^3-H}$.

A nonseparating disk has distance $1$ if and only if it is primitive, since
both conditions are equivalent to the condition that cutting $H$ along the
disk produces an unknotted solid torus (see \cite[Section
4]{JohnsonBridgeNumber}). Therefore the tunnel of the trivial knot is the
only tunnel of distance~$1$. A simple or semisimple tunnel has
distance~$2$, since it is disjoint from a primitive disk. There are,
however, regular tunnels of distance~$2$. In section~\ref{sec:torus_knots}
we will see that the regular tunnels of torus knots all have distance~$2$.

Recall that if $\Sigma=(\overline{H-\Nbd(K_\tau)},\overline{S^3-H})$ is a
Heegaard splitting of the complement of $K_\tau$, then the \textit{(Hempel)
distance $\dist(\Sigma)$} is the minimal distance in the curve complex of
$\partial H$ between the boundary of a disk in $\overline{H-\Nbd(K_\tau)}$
and the boundary of a disk in $\overline{S^3-H}$ (where the disks may be
separating). Clearly, $\dist(\Sigma)\leq \dist(\tau)$. On the other hand,
Johnson \cite[Lemma 11]{JohnsonBridgeNumber} proved that
\begin{lemma}[Johnson] $\dist(\tau)\leq \dist(\Sigma)+1$.
\label{lem:JohnsonDistance}
\end{lemma}

M. Scharlemann and M. Tomova~\cite{Scharlemann-Tomova} proved the following
stability result:
\begin{theorem}[Scharlemann-Tomova]
Genus-$g$ Heegaard splittings of distance more than $2g$ are isotopic.
\label{thm:Scharlemann-Tomova}
\end{theorem}
\noindent Using lemma~\ref{lem:JohnsonDistance} and
theorem~\ref{thm:Scharlemann-Tomova}, Johnson~\cite[Corollary
13]{JohnsonBridgeNumber} deduced the following:
\begin{theorem}[Johnson] If $\tau$ is a tunnel of a tunnel number $1$ knot
$K_\tau$ and $\dist(\tau)>5$, then $\tau$ is the unique tunnel of $K_\tau$.
\label{thm:uniquetunnel}
\end{theorem}

Theorem~15.2 of~\cite{CM}, an immediate consequence of the Parameterization
Theorem, determines all orientation-reversing self-equivalences of tunnels:
\begin{theorem} Let $\tau$ be a tunnel of a tunnel number~$1$ knot or link.
Suppose that $\tau$ is equivalent to itself by an orientation-reversing
equivalence. Then $\tau$ is the tunnel of the trivial knot, the trivial
link, or the Hopf link.
\label{thm:HopfLink}
\end{theorem}
Combining theorems~\ref{thm:uniquetunnel} and~\ref{thm:HopfLink} gives the
following:
\begin{corollary}
If $\tau$ is a tunnel of a tunnel number $1$ knot $K_\tau$ and
$\dist(\tau)>5$, then $K_\tau$ is not amphichiral.
\label{coro:amphichiral}
\end{corollary}
\noindent For theorem~\ref{thm:HopfLink} shows that an
orientation-reversing equivalence from $K_\tau$ to $K_\tau$ would produce a
second tunnel for~$K_\tau$.

Distance also has implications for hyperbolicity.
\begin{theorem} If $K_\tau$ is a torus knot or a satellite knot, then 
$\dist(\tau)\leq 2$.
Consequently, if $\dist(\tau)\geq 3$, then $K_\tau$ is hyperbolic.\par
\label{thm:Morimoto-Sakuma}
\end{theorem}
\begin{proof} We have already mentioned the fact, verified in
section~\ref{sec:torus_knots} below, that the regular tunnels of torus
knots have distance~$2$. The other tunnels of torus knots are simple or
semisimple, so also have distance~$2$. K. Morimoto and M. Sakuma
\cite{Morimoto-Sakuma} found all tunnels of tunnel number~$1$ satellite
knots, showing in particular that they are semisimple.
\end{proof}

The \textit{depth} of $\tau$ is the simplicial distance $\depth(\tau)$ in
the $1$-skeleton of $\D(H)/\G$ from $\tau$ to the primitive vertex $\pi_0$.
The inequality
\[\dist(\tau)-1\leq \depth(\tau)\]
mentioned in the introduction is immediate from the definitions.  Also from
the definitions, $\tau$ is primitive if and only if $\depth(\tau)=0$, is
simple or semisimple if and only if $\depth(\tau)=1$, and is regular if and
only if $\depth(\tau)\geq 2$. As already mentioned, in
section~\ref{sec:torus_knots} we will see a sequence of torus knot tunnels
of distance~$2$ with depths that grow arbitrarily large.

\section{Giant steps} 
\label{sec:GST_moves}

\begin{definition} Let $\tau$ and $\tau'$ be tunnels tunnel number~$1$ knot.
We say that $\tau'$ is obtained from $\tau$ by a \textit{giant step} if
$\tau$ and $\tau'$ are the endpoints of a $1$-simplex of
$\D(H)/\G$. Equivalently, $\tau$ and $\tau'$ can be represented by disjoint
disks in~$H$.
\end{definition}
It is clear that the depth of a tunnel is the minimum number of giant steps
needed to transform the tunnel to the tunnel of the trivial knot, or vice
versa.

In \cite{GST}, Goda, Scharlemann, and Thompson gave a geometric definition
of giant steps (this is one reason for our the selection of the name Giant
STeps), as follows. Let $\tau$ be a nonseparating disk in $H$, and let $K$
be a simple closed curve in $\partial H$ that intersects $\tau$
transversely in one point. Let $N$ be a regular neighborhood in $H$ of
$K\cup \gamma$. Then the frontier of $N$ separates $H$ into two solid tori,
one a regular neighborhood of $K$, so $K$ is a tunnel number~$1$ knot.

In the previous construction, the meridian disk $\tau'$ of the solid torus
that does not contain $K\cup\tau$ is the \textit{unique} nonseparating disk
$\tau'$ in $H$ that is disjoint from $K\cup \tau$, and $\tau'$ is a tunnel
of $K$. That is, the construction produces a specific tunnel of the
resulting knot $K$. A giant step as we have defined it simply amounts to
choosing the $\tau'$ first; $K$ is then determined up to isotopy in $H$ and
in $S^3$, although not up to isotopy in~$\partial H$.

Since the complex $\D(H)/\G$ is connected, we have the following, which is
part of Proposition~1.11 of~\cite{GST}.
\begin{proposition}
Let $\tau$ be a tunnel of a tunnel number $1$ knot. Then there is a
sequence of giant steps that starts with the tunnel of the trivial knot and
ends with~$\tau$.\par
\label{prop:GST_sequence_exists}
\end{proposition}

\section{Minimal sequences of giant steps}
\label{sec:GST_move_seqs}

In this section, we will calculate the number of minimal length sequences
of giant steps that start from $\pi_0$, the tunnel of the trivial knot, and
end with a given tunnel~$\tau$. First, we will observe that any such
minimal sequence corresponds to a minimal length simplicial path in the
$1$-skeleton of a neighborhood of the principal path of $\tau$. An
elementary counting method then gives an algorithm to calculate the number
of such paths.

By a \textit{path} (between two vertices) in $\D(H)/\G$, we mean a
simplicial path in the $1$-skeleton of $\D(H)/\G$, passing through a
sequence of vertices that are images of vertices of $\D(H)$ (i.~e.~vertices
that represent tunnels).  We describe such a path simply by listing the
vertices through which it passes. From section~\ref{sec:GST_moves}, we know
that the minimal sequences of giant steps from the trivial tunnel to a given
tunnel correspond exactly to the minimal-length paths in $\D(H)/\G$ from
$\pi_0$ to the tunnel vertex. We will only be interested in minimal-length
paths.

\begin{definition}
Let $\tau$ be a nontrivial tunnel. Define the \textit{corridor} of $\tau$,
$C(\tau)$, as follows. Write the vertices of the principal path of $\tau$
as $\theta_0$, $\mu_0$, $\mu_0\cup \tau_0$, $\mu_1\,$, $\mu_1\cup
\tau_1,\ldots\,$, $\mu_n\cup \tau_n$, where $\tau=\tau_n$. Then $C(\tau)$
is the union of the $2$-simplices whose barycenters are the $\mu_i\cup
\tau_i$ for $0\leq i\leq n$ (where $\mu_0\cup \tau_0$ is regarded as the
barycenter of the $2$-simplex spanned by $\pi_0$, $\mu_0$, and $\tau_0$).
\label{def:corridor}
\end{definition}
\begin{figure}
\begin{center}
\includegraphics[width=66 ex]{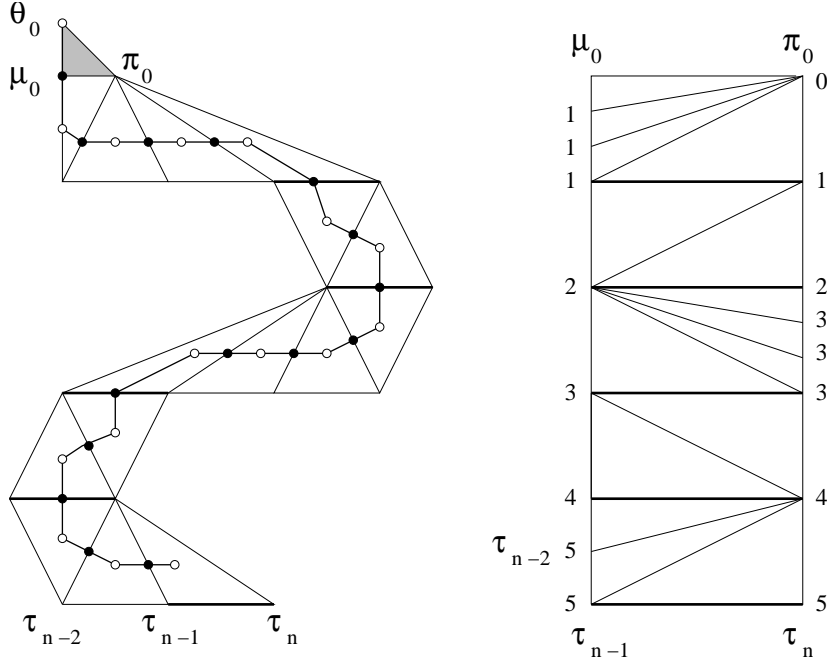}
\caption{The principal path and corridor $C(\tau_n)$ of a tunnel
$\tau_n$. The $\nabla$-edges are emphasized. In the picture of
$C(\tau_n)$ on the right, the depths of the tunnels are labeled.}
\label{fig:corridor}
\end{center}
\end{figure}

When $\tau$ is a simple tunnel, $C(\tau)$ is the triangle
$\langle\pi_0,\mu_0,\tau_0\rangle$. Otherwise, it can be viewed as a
rectangular or trapezoidal strip with end edges $\langle
\pi_0,\mu_0\rangle$ and $\langle \tau_j,\tau_n\rangle$ for some $j$, as in
the figure on the right in figure~\ref{fig:corridor}.

\begin{lemma} Let $\tau$ be a tunnel, and let
$\sigma_0$, $\sigma_1,\ldots\,$, $\sigma_n$ be a path in $\D(H)/\G$ of
minimal length among the paths connecting $\sigma_0$ to $\sigma_n$.  If
$\sigma_0$ and $\sigma_n$ lie in $C(\tau)$, then each $\sigma_i$ lies
in~$C(\tau)$.\par
\label{lem:stay_in_corridor}
\end{lemma}

\begin{proof} 
If the lemma is false, then there exist $i$ and $j$ with $0<i<i+1<j<n$ for
which $\sigma_i$ and $\sigma_j$ lie in the frontier in $\D(H)/\G$ of
$C(\tau)$, but $\sigma_k$ does not lie in $C(\tau)$ for any $m$ with
$i<k<j$.  Let $\sigma_i'$ and $\sigma_i''$ be the vertices of $\D(H)/\G$
adjacent to $\sigma_i$ on the frontier of $C(\tau)$. The union of
$1$-simplices $\langle \sigma_i',\sigma_i\rangle$ and $\langle
\sigma_i,\sigma_i''\rangle$ separates $\D(H)/\G$, indeed every $1$-simplex
of $\D(H)/\G$ separates. Therefore $\sigma_j$ must equal one of $\sigma_i$,
$\sigma_i'$, or $\sigma_i''$. But this implies that the original path did
not have minimal length.
\end{proof}

In the special case that $\tau$ is of depth $1$, $\tau$ lies in the link in
$\D(H)/\G$ of $\pi_0$, and there is clearly a unique path of length~$1$
from $\tau$ to $\pi_0$. From now on, we assume that $\tau$ has depth at
least~$2$.

Now, regard $C(\tau)$ as in the diagram on the right in
figure~\ref{fig:corridor}, with the edge $\langle \mu_0, \pi_0\rangle$ on
top, and with $\tau$ as one of the endpoints of the bottom edge.  In the
triangulation of $C(\tau)$, a \textit{$\nabla$-edge of depth $n$} is an
edge whose endpoints have depth $n$ and lie on different sides of
$C(\tau)$, and for which all vertices lying below its endpoints on either
side have depth greater than $n$. In figure~\ref{fig:corridor}, the
$\nabla$-edges are highlighted.

Since the endpoints of any edge of $C(\tau)$ can have depths that differ by
at most $1$, there exists a unique $\nabla$-edge $\nabla(i)$ in $C(\tau)$
of depth $i$ for each $i$ with $1\leq n<\depth(\tau)$ (and there may be one
of depth~$n$).

The name $\nabla$-edges arises from the fact that these edges are the tops
of $2$-simplices of the corridor that appear as $\nabla$'s when the
corridor is drawn with depth corresponding to the vertical coordinate, as
in the diagram on the left in figure~\ref{fig:corridor}. Every nonprimitive
$2$-simplex of $\D(H)/\G$ has two vertices of the same depth and a third of
depth either $1$ larger or $1$ smaller than that depth; for a ``$\nabla$''
$2$-simplex that depth is $1$ larger, while it is $1$ smaller for a
``$\Delta$'' $2$-simplex.

Denote the left and right endpoints of $\nabla(i)$ by
$\partial_L(\nabla(i))$ and $\partial_R(\nabla(i))$ respectively.
\begin{lemma} Let $\nabla(i-1)$ and $\nabla(i)$ be successive
$\nabla$-edges. Then at least one of the pairs
$\{\partial_L(\nabla(i-1)),\partial_L(\nabla(i))\}$ and 
$\{\partial_R(\nabla(i-1)),\partial_R(\nabla(i))\}$ are
the endpoints of an edge that lies in a side of~$C(\tau)$.
\label{lem:side_edge}
\end{lemma}
\begin{proof}
For each endpoint of $\nabla(i)$, select a path of length $i$ from the
endpoint to $\pi_0$.  By lemma~\ref{lem:stay_in_corridor}, these paths lie
in $C(\tau)$. In particular, each of their first edges connects an endpoint of
$\nabla(i)$ to an endpoint of $\nabla(i-1)$. At most one of these first
edges can be diagonal, so at least one lies in a side.
\end{proof}

\begin{figure}
\begin{center}
\includegraphics[width=\textwidth]{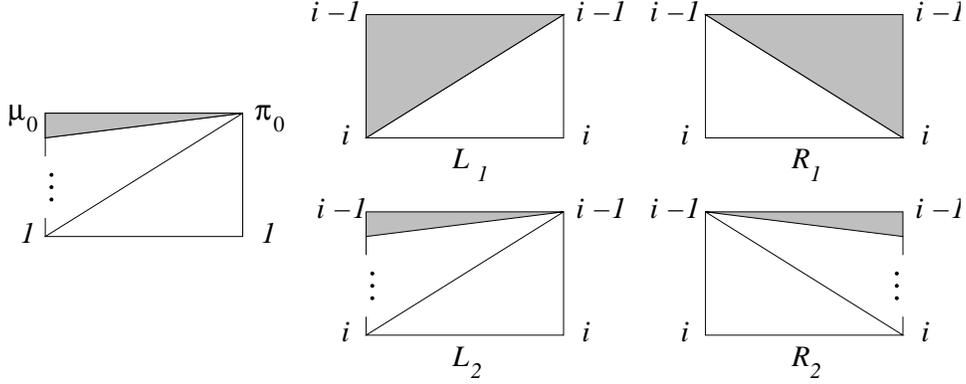}
\caption{The configuration above $\nabla(1)$, and
the four possible configurations between two $\nabla$-edges.
In $L_1$ and $R_1$ there is only one $2$-simplex above the diagonal edge,
while in $L_2$ and $R_2$ there are two or more. In the configuration
above $\nabla(1)$, there may be only one $2$-simplex above the diagonal.
The shaded $2$-simplices are $\nabla$ $2$-simplices.
The letter $L$ (respectively, $R$) signifies that the $\nabla$ $2$-simplex
has an edge on the left side (respectively, right side) of the corridor.}
\label{fig:configurations}
\end{center}
\end{figure}

Lemma~\ref{lem:side_edge} shows that the triangulation of the portion of
$C$ between $e$ and $f$ must have one of the four configurations $L_1$,
$R_1$, $L_2$, or $R_2$ shown in figure~\ref{fig:configurations}. The
portion of $C(\tau)$ above $\nabla(1)$ must be as in the leftmost diagram
in figure~\ref{fig:configurations}, where there may be only one $2$-simplex
above the diagonal.

Now, we show how to calculate the number of minimal paths from $\tau$ to
$\pi_0$. Denote by $\lambda_i$ the number of paths in $C(\tau)$ of length
$\depth(\tau)-i$ from $\tau$ to the left endpoint of $\nabla(i)$, and by
$\rho_i$ the number to its right endpoint. From
figure~\ref{fig:configurations}, we see that the total number of paths from
$\tau$ to $\pi_0$ is $\lambda_1+\rho_1$.

Let $\nabla(n)$ be the $\nabla$-edge in $C(\tau)$ of maximum depth.  If
$\tau$ is the left (or right) endpoint of $\nabla(n)$, then
$\begin{pmatrix} \lambda_n \\ \rho_n \end{pmatrix}$ is $\begin{pmatrix} 1 \\ 0
\end{pmatrix}$ (or $\begin{pmatrix} 0\\ 1 \end{pmatrix}$).  If $\tau$ and
the endpoints of $\nabla(n)$ span a $2$-simplex, as in the case of the
tunnel $\tau_{n-2}$ in figure~\ref{fig:corridor}, then $\begin{pmatrix}
\lambda_n \\ \rho_n
\end{pmatrix}$ is $\begin{pmatrix} 1 \\ 1 \end{pmatrix}$.
Otherwise, there
is exactly one edge from $\sigma$ to an endpoint of $\nabla(n)$, and
$\begin{pmatrix} \lambda_n \\ \rho_n
\end{pmatrix}$ is $\begin{pmatrix} 1 \\ 0 \end{pmatrix}$ or
$\begin{pmatrix} 0\\ 1 \end{pmatrix}$ according to the endpoint.

For each $2\leq i\leq n$, let $C_i$ be $L_1$, $R_1$, $L_2$, or $R_2$
according to which of the four configurations in
figure~\ref{fig:configurations} describes the triangulation of $C$ between
$\nabla(i-1)$ and $\nabla(i)$.  For $2\leq i\leq n-1$, put $M_i$ equal to
the matrix given in the following table, according to the value of $C_i$:
\medskip

\begin{small}
\setlength{\fboxsep}{0pt}
\setlength{\fboxrule}{0.5pt}  
\begin{center}
\fbox{%
\begin{tabular}{c|c|c|c|c}
$C_i$&$L_1$&$R_1$&$L_2$&$R_2$\\
\hline
$M_i$&$\begin{pmatrix}1 & 0\\ 1 & 1 \end{pmatrix}$
&$\begin{pmatrix}1 & 1\\ 0 & 1 \end{pmatrix}$
&$\begin{pmatrix}0 & 0\\ 1 & 1 \end{pmatrix}$
&$\begin{pmatrix}1 & 1\\ 0 & 0 \end{pmatrix}$\\
\end{tabular}}
\end{center}
\end{small}
\medskip

Observe that
\[ M_i\begin{pmatrix} \lambda_i \\ \rho_i\end{pmatrix} =
\begin{pmatrix} \lambda_{i-1} \\ \rho_{i-1}\end{pmatrix}\ .\]
Therefore we have
\[ \begin{pmatrix} \lambda_1 \\ \rho_1\end{pmatrix} =
M_2M_3\cdots M_n\begin{pmatrix} \lambda_n \\ \rho_n\end{pmatrix}\ ,\]
or alternatively, putting $M_1=\begin{pmatrix}1 & 1\end{pmatrix}$, the number
of minimal paths from $\tau$ to $\pi_0$ is the entry of the
$1\times 1$-matrix $M_1M_2\cdots M_n \begin{pmatrix} \lambda_n \\ \rho_n\end{pmatrix}$.

The algorithm just described is easily implemented
computationally~\cite{slopes}. The input is the binary string $s_2s_3\cdots
s_n$ of parameters from the Parameterization Theorem, which completely
determine the structure of $C(\tau)$. The input is broken into blocks
having one of the four forms $10$, $11$, $100+$, and $110+$, where $0+$
indicates a nonempty string of $0$'s. A $10$-block, for example, produces a
configuration $A_1$ when the principal path is moving in from right to
left, and a configuration $B_1$ when it is moving from left to right, and
so on. The blocks $10$ and $100+$ reverse the direction, and the others do
not. A leftover $1$ at the end of the input string indicates that the
bottom triangle in $C(\tau)$ is a $\nabla$~$2$-simplex, so that
$\begin{pmatrix} \lambda_n \\ \rho_n
\end{pmatrix}=
\begin{pmatrix} 1 \\ 1\end{pmatrix}$, and in the other cases
$\begin{pmatrix} \lambda_n \\ \rho_n \end{pmatrix}$ is worked out from the
final block (and the direction of travel of the principal path at that
point).

For the example in figure~\ref{fig:corridor}, the input string is
$0011100011100$ and the output of the program is:
\medskip

\begin{ttfamily}
Depth> gst( '0011100011100', verbose=True )

The intermediate configurations are L1, R2, R1.

The transformation matrices are:

\ \ \  [ [ 1, 0 ], [ 1, 1 ] ]

\ \ \   [ [ 1, 1 ], [ 0, 0 ] ]

\ \ \   [ [ 1, 1 ], [ 0, 1 ] ]

and their product is [ [ 1, 2 ], [ 1, 2 ] ].

The final block has configuration L2.

This tunnel has 4 minimal giant step constructions.
\end{ttfamily}
\medskip

Among the interesting examples are the tunnels whose parameter sequences
are the following:
\begin{enumerate}
\item $s_2s_3\cdots s_n=100100\cdots 100$. The configuration
sequence alternates as $R2$, $L2$, $R2$, $L2\ldots\,$, and there is a
unique minimal giant step sequence.
\item $s_2s_3\cdots s_{2n+1}=1010\cdots 10$. The configuration
sequence alternates as $R1$, $L1$, $R1$, $L1\ldots\,$, and the number
of minimal giant step sequences is $a_n$, the term in the Fibonacci sequence
$(a_0,a_1,a_2,\ldots)=(1,1,2,3,5,\ldots)$.
\item $s_2s_3\cdots s_{2n+1}=111\cdots 1$, an even number of $1$'s.
The configuration sequence is
$L1$, $L1\ldots\,$, and there is a unique minimal giant step sequence.
\item $s_2s_3\cdots s_{2n}=111\cdots 1$, an odd number of $1$'s. The
configuration sequence is $L1$, $L1\ldots\,$, but there is only a single
$\nabla$ $2$-simplex below $\nabla(n)$, and there are $n+1$ minimal giant
step sequences.
\end{enumerate}
Note that examples of the last two types are obtained from each other by a
single additional cabling construction, even though the numbers of minimal
giant step constructions differ by arbitrarily large amounts.

\begin{remark} The algorithm shows that tunnels with a unique minimal
giant step sequence are sparse. For instance, the product $R_1L_1R_1L_1$ is
a matrix with all entries greater than~$1$, so whenever this appears as any
block of four terms in the product $M_1M_2\cdots M_n$ that occurs in the
algorithm, there must be more than one minimal giant step
sequence. Products $M_1M_2\cdots M_n$ in which this block occurs are
generic in any reasonable sense.
\end{remark}

\section{Depth and bridge number}
\label{sec:bridge_number_growth}

Some deep geometric results of Goda, Scharlemann, and Thompson allow us to
obtain information about the bridge numbers of the knots $K_\tau$. They
show that once one leaves the semisimple region, the bridge number grows at
least exponentially with the depth (in fact, it grows rapidly with the
number of cablings needed to produce the tunnel).

\begin{figure}
\begin{center}
\includegraphics[height=42 ex]{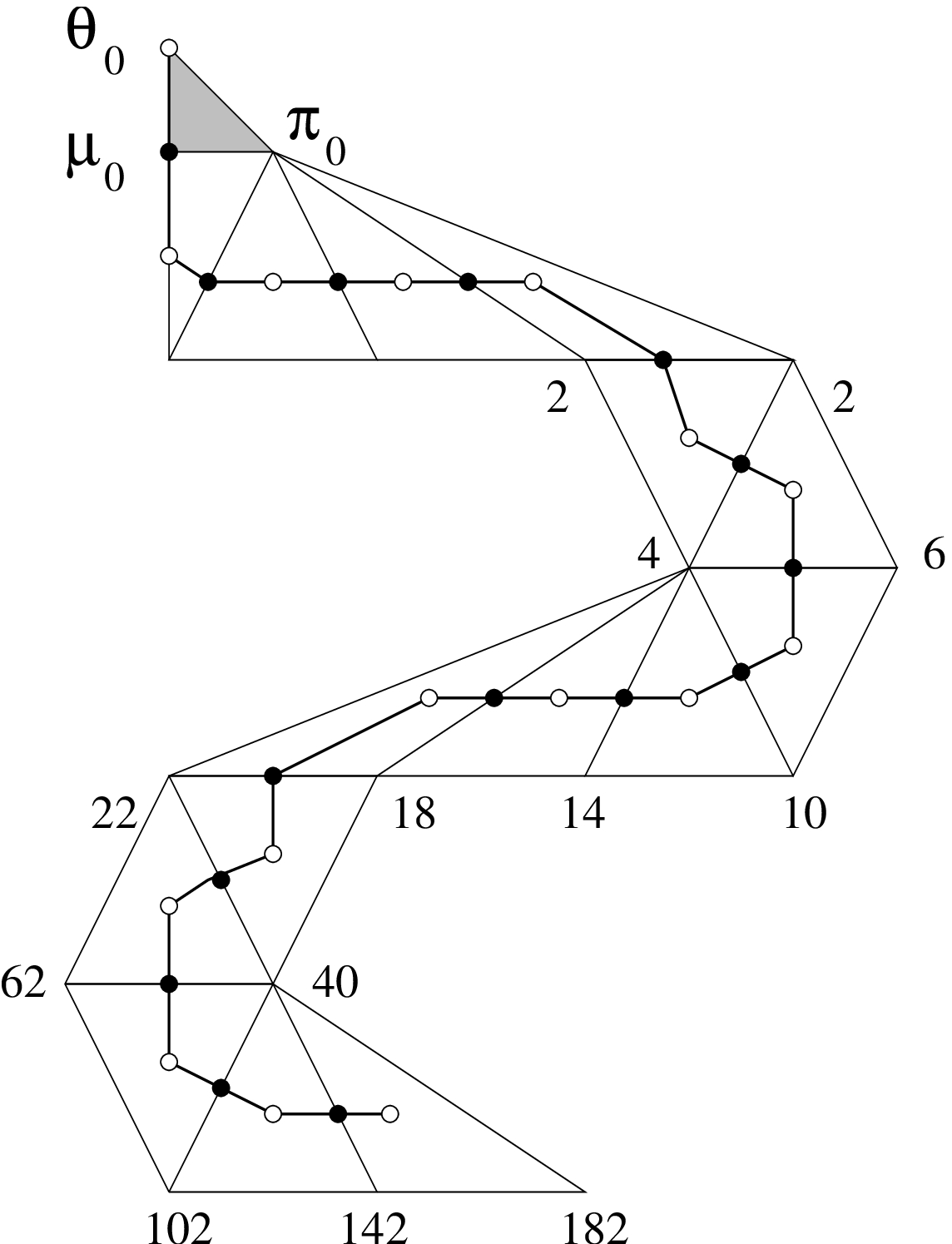}
\caption{}
\label{fig:bridge_numbers}
\end{center}
\end{figure}

Our most precise result on bridge numbers is
theorem~\ref{thm:bridge_numbers}, which gives a lower bound for the bridge
number of a tunnel number~$1$ knot in terms of the principal path of any of
its regular tunnels. Its statement is a bit uninviting, but it says
something very simple. Figure~\ref{fig:bridge_numbers} illustrates how
theorem~\ref{thm:bridge_numbers} bounds the bridge numbers of the tunnels
along the principal path of the example of figure~\ref{fig:corridor},
assuming that the last two tunnels at depth~$1$ were tunnels of $2$-bridge
knots.
\begin{theorem} Let $\tau_n$ be
a regular tunnel with principal path $\theta_0$, $\mu_0$, $\mu_0\cup
\{\tau_0\}$, $\mu_1,\ldots\,$, $\mu_n$, $\mu_n\cup \{\tau_n\}$,
$\tau_n$. In the principal path of $\tau$, let $\tau_m$ be the first tunnel
of depth $2$, with principal vertex $\{\tau_{m-2},\tau_{m-1},\tau_m\}$. Put
$c_{m-2}=\br(K_{\tau_{m-2}})$ and $c_{m-1}=\br(K_{\tau_{m-1}})$. To the
vertices $\tau_m$, $\tau_{m+1},\ldots\,$, $\tau_n$, assign values $c_k$
inductively by the rule $c_k=c_i+c_j$, where the principal vertex of
$\tau_k$ is $\{\tau_i,\tau_j,\tau_k\}$. Then the bridge number of
$K_{\tau_n}$ is at least~$c_n$.
\label{thm:bridge_numbers}
\end{theorem}

Before proving theorem~\ref{thm:bridge_numbers}, we isolate the step that
uses the results of Goda, Scharlemann, and Thompson~\cite{GST}
and Scharlemann and Thompson~\cite{Scharlemann-Thompson}.
For its statement and proof, we remind the reader
that the term \textit{principal meridian pair} was defined near the
end of section~\ref{sec:disk_complex}.
\begin{lemma} Let $\{\lambda,\rho\}$ be the principal meridian pair of a
tunnel $\tau$, and let $\theta$ be the $\theta$-curve associated 
to the principal vertex $\{\lambda,\rho, \tau\}$ of $\tau$. Write $T$ for
the arc dual to $\tau$, and $L$ and $R$ for the other two arcs of $\theta$
that are dual to $\lambda$ and $\rho$, so that $K_{\tau} = L \cup R$,
$K_\lambda=R \cup T$, and $K_\rho=L \cup T$. Then
there a minimal bridge position of $K_\tau$ in either:
\begin{enumerate}
\item[(i)] $T$ is slid to an arc in a level sphere, and $T$ connects two
bridges of $K_\tau$. In the $n$-strand trivial tangle above the level
sphere, the arcs are parallel to a collection of disjoint arcs in the level
sphere, which meet $T$ only in its endpoints. Moreover, $K_\tau\cup T$ is
isotopic to the original $\theta$. Or,
\item[(ii)] $T$ is slid to an eyeglass in a level sphere. The endpoints of
$T$ can be slid slightly apart, moving $T$ out of the level sphere,
producing $K_\tau\cup T$ isotopic to the original $\theta$, and showing
that one of $K_\lambda$ or $K_\rho$ is a trivial knot, and consequently
$\tau$ is simple or semisimple.
\end{enumerate}
\label{lem:GSTlemma}
\end{lemma}
\begin{proof}
By Theorem~1.8 of \cite{GST}, we may move $K_\tau \cup T$, possibly using
slide moves of $T$ as well as isotopy, so that $K_\tau$ is in minimal
bridge position and $T$ either lies on a level sphere and connects two
bridges of $K_\tau$, or $T$ is slid to an ``eyeglass''. Since the leveling
process involves sliding the tunnel arc $T$, there is \textit{a priori} no
reason for the resulting $\theta$-curve to be isotopic to the original
$\theta$. But Corollary~3.4 and Theorem~3.5 (combined with Lemma~2.9) of
\cite{Scharlemann-Thompson} show that in~(i) and~(ii), the dual disks to
the other two arcs of the $\theta$-curve are the principal meridian pair of
$\tau$, that is, $\lambda$ and $\rho$, so the resulting $\theta$-curve is
still~$\theta$. Finally, the description of the trivial tangle above the
level sphere in case~(i) is from Theorem~6.1 of~\cite{GST}.
\end{proof}

Theorem~\ref{thm:bridge_numbers} follows immediately by
application of the next lemma, which will also be used in the proof of
corollary~\ref{coro:bridge_numbers}.
\begin{lemma} Let $\tau$ be a tunnel of a nontrivial knot, and
let $\{\lambda,\rho\}$ be the principal meridian pair of $\tau$.  Then
$\br(K_\tau)\geq \br(K_\lambda)+\br(K_\rho)-1$. If $\tau$ is regular,
then $\br(K_\tau)\geq \br(K_\lambda)+\br(K_\rho)$.
\label{lem:bridge_number_inequality}
\end{lemma}
\begin{proof} Level a tunnel arc as in lemma~\ref{lem:GSTlemma}. If the
tunnel arc is slid to an eyeglass, pulling the endpoints slightly apart
produces $\theta$ and shows that one of $K_\lambda$ or $K_\rho$ is
trivial. The other is in bridge position with the same number of bridges as
$K_\tau$, so $\br(K_\tau)\geq \br(K_\lambda)+\br(K_\rho)-1$. If $\tau$ is
regular, then the eyeglass configuration cannot not occur, and we see that
$\br(K_\tau)\geq \br(K_\lambda)+\br(K_\rho)$.
\end{proof}

It is straightforward to implement the iteration of
theorem~\ref{thm:bridge_numbers} computationally~\cite{slopes}. The only
information needed for input is the sequence of parameters $s_2,\ldots\,$,
$s_n$ of the Parameterization~Theorem and the values $c_{m-2}$ and
$c_{m-1}$, called $c_2$ and $c_3$ in the sample output shown here:
\medskip

\begin{ttfamily}
Depth> bridge( '0011100011100', c2=2, c3=2, verbose=True )

The minimum bridge number of K-tau is 182.

The iteration sequence is:

\ \ \ 2, 2, 4, 6, 10, 14, 18, 22, 40, 62, 102, 142, 182
\end{ttfamily}
\medskip

Theorem~\ref{thm:bridge_numbers} gives us a general lower bound for bridge
number as a function of depth:
\begin{corollary} Let $\tau$ be a regular tunnel of depth $d$, and in
the principal path of $\tau$, let $\tau_m$ be the first tunnel of depth
$2$, with principal vertex $\{\tau_{m-2},\tau_{m-1},\tau_m\}$. Put
$b_2=\br(K_{\tau_{m-2}})$ and $b_3=\br(K_{\tau_{m-1}})$. For $n\geq 2$ let
$b_j$ be given by the recursion
\begin{gather*} b_{2n}=b_{2n-1}+b_{2n-2}\\
b_{2n+1}=b_{2n}+b_{2n-2}
\end{gather*}
Then $\br(K_\tau)\geq b_{2d}$.
\label{coro:bridge_numbers}
\end{corollary}
\begin{figure}
\begin{center}
\includegraphics[height=35 ex]{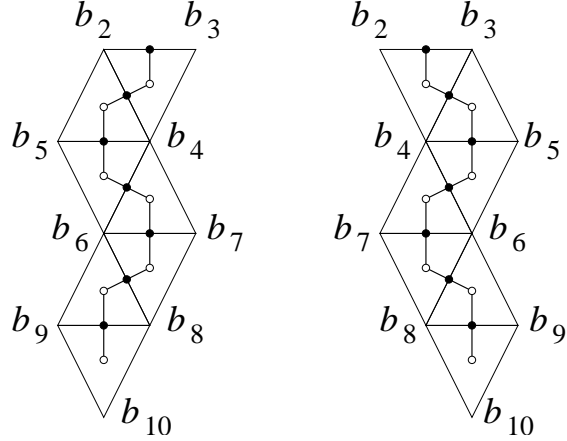}
\caption{The path on the left is of cheapest descent. The path on the right
is of cheapest descent if $b_2=b_3$.}
\label{fig:descent}
\end{center}
\end{figure}
\begin{proof} By an application of
lemma~\ref{lem:bridge_number_inequality}, we have $b_2\leq b_3$.  The left
diagram in figure~\ref{fig:descent} shows a ``path of cheapest descent''
starting from the vertex $\{\tau_{m-2},\tau_{m-1}\}$ (the path in the right
diagram is also a path of cheapest descent, provided that $b_2=b_3$).  Any
principal path having more than two tunnels at a given depth will produce
an even larger bridge number, as will any principal path that emerges in
the more costly direction out of a $\nabla$ $2$-simplex. Applying
theorem~\ref{thm:bridge_numbers} to a path of cheapest descent gives the
recursion of corollary~\ref{coro:bridge_numbers}, hence a lower bound
for~$\br(K_\tau)$.
\end{proof}

We can now prove one of our main results.
\begin{minimumBridgeNumberTheorem} 
For $d\geq 1$, the minimum bridge number of a knot having a tunnel of
depth~$d$ is $a_d$, where $a_1=2$, $a_2=4$, and $a_d=2a_{d-1}+a_{d-2}$ for
$d\geq 3$.
\label{thm:minbridgenum}
\end{minimumBridgeNumberTheorem}
\begin{proof} Taking $b_2=b_3=2$ in corollary~\ref{coro:bridge_numbers}
gives a $b_{2d}$ which is a general lower bound for the bridge number of a
tunnel at depth~$d$, and a little bit of algebra shows that $b_{2d}=a_d$
for the recursion in theorem~\ref{thm:minbridgenum}. It remains to show the
existence of a $\tau$ of depth $d$ for which $\br(K_\tau)=a_{2d}$.

We begin with a tunnel $\rho$ which is a semisimple tunnel of a $2$-bridge
knot obtained from the trivial tunnel by two cablings. Details of the
construction of such a tunnel are given in section~17 of~\cite{CM}.  The
principal path of $\rho$ is a portion of the path shown in the leftmost
diagram of figure~\ref{fig:bridgenum4}, and the principal pair of $\rho$ is
$\{\pi_0,\lambda\}$. A tunnel $\tau$ is constructed using a cabling as
indicated in figure~\ref{fig:bridgenum4}, producing a $4$-bridge knot
$K_\tau$. The number of twists in the two strands on the right is variable,
it must simply be chosen so that $K_\tau$ is a knot rather than a link. The
slope of this cabling is an odd integer, since $K_\tau$ meets the replaced
disk $\pi_0$ in only two points.
\begin{figure}
\begin{center}
\includegraphics[width=\textwidth]{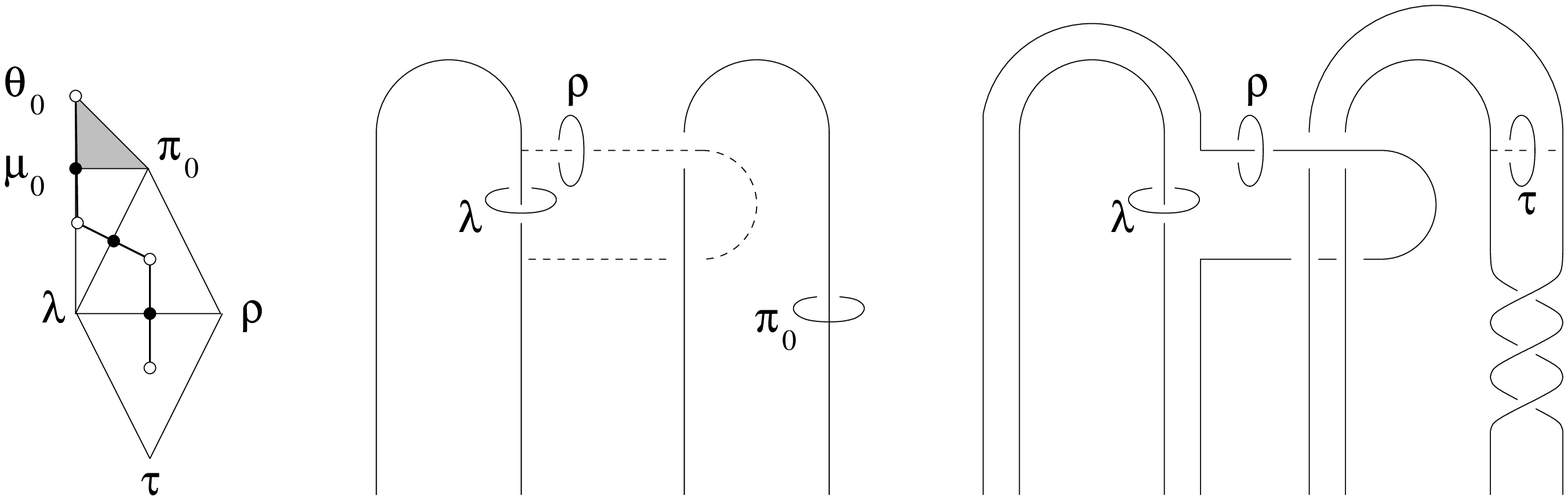}
\caption{}
\label{fig:bridgenum4}
\end{center}
\end{figure}

\begin{figure}
\begin{center}
\includegraphics[width=\textwidth]{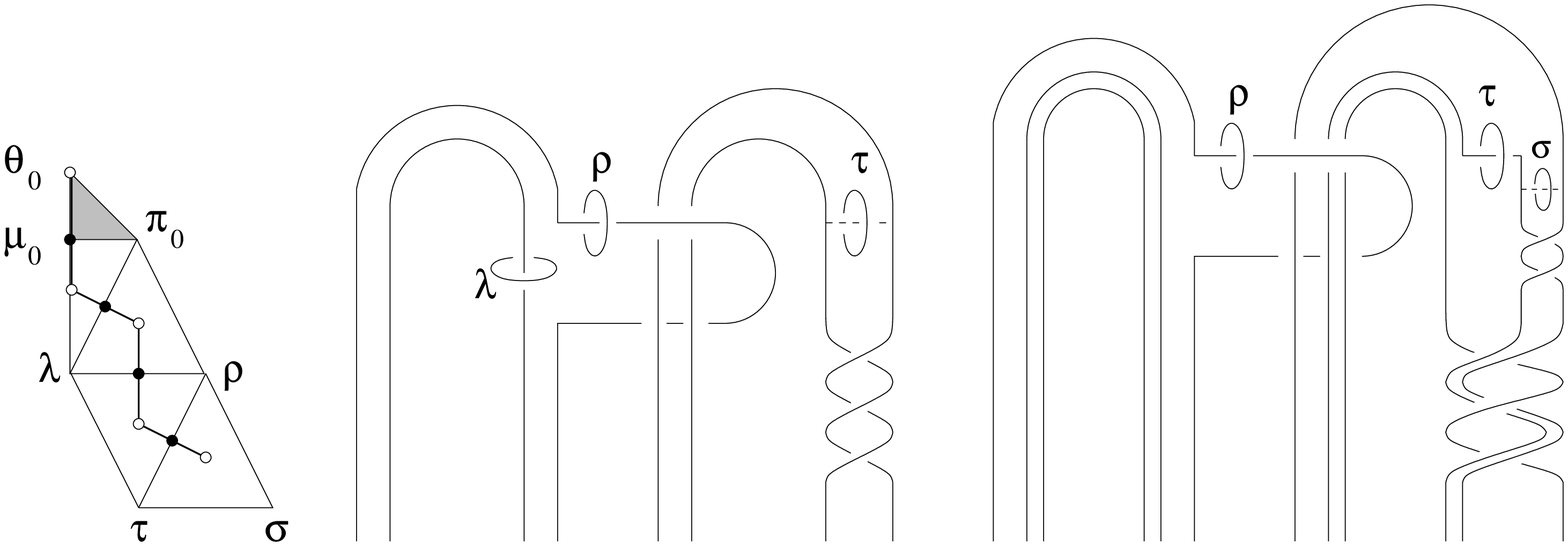}
\caption{}
\label{fig:bridgenum6}
\end{center}
\end{figure}
Figure~\ref{fig:bridgenum6} shows how to continue the construction. The
next tunnel $\sigma$ must be at depth~$2$, and at this stage (since
$b_2=b_3$) we can retain either $\lambda$ or $\rho$; we have chosen to
retain $\rho$ in the example of figure~\ref{fig:bridgenum6}. Again, the
cabling has odd integer slope. The resulting tunnel $\sigma$ has
$\br(K_\sigma)=6$. Later repetitions of the construction resemble that of
figure~\ref{fig:bridgenum6}, but the analogues of $K_{\lambda}$ and
$K_{\rho}$ will not have the same bridge number, and one must retain
whichever disk has corresponding knot of smaller bridge number. A pattern
as in figure~\ref{fig:descent} will be produced.
\end{proof}

\section{Upper bounds for bridge number}
\label{sec:bridgenum_conj}

\begin{figure}
\begin{center}
\includegraphics[width=\textwidth]{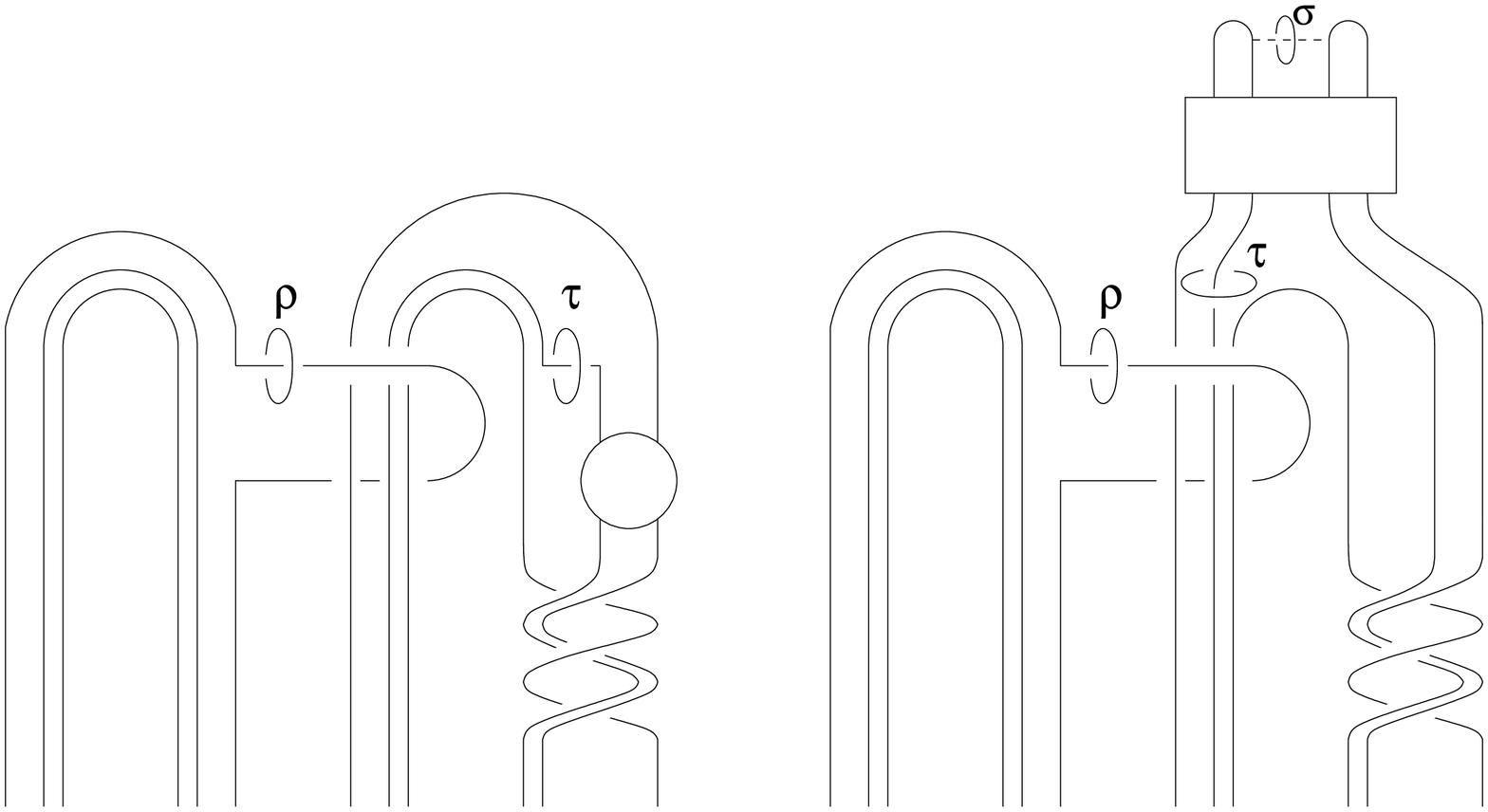}
\caption{}
\label{fig:minraise}
\end{center}
\end{figure}

The construction used to prove the Minimum Bridge Number Theorem adapts to
show that in general, cabling operations can be carried out efficiently
from the viewpoint of bridge number. The basic idea is shown in
figure~\ref{fig:minraise}.  Its left diagram is like the right diagram of
figure~\ref{fig:bridgenum6}, except that a cabling of some arbitrary slope
has been performed on $K_\tau$ to produce $K_\sigma$; the rational tangle
in $K_\sigma$ created by the cabling is inside a ball represented by the
circle in the first diagram. We can reposition $K_\sigma$ as in the second
diagram of figure~\ref{fig:minraise}, by ``moving the ball up to engulf
infinity,'' in such a way that the rectangle in the second diagram contains
a $4$-strand braid.  The new tunnel $\sigma$ is still level so the
construction can be repeated.

To understand the effect of this construction on bridge numbers, we
introduce some special terminology. 
\begin{definition} Let $\{\lambda,\rho,\tau\}$ be the principal vertex of a
tunnel $\tau$. Its dual $\theta$-curve has the form $K_\tau\cup \alpha$,
where $\alpha$ is a tunnel arc representing the tunnel $\tau$, and contains
$K_\lambda$ and $K_\rho$ as subsets. Assume that $K_\tau\cup \alpha$ is
positioned so that $K_\tau$ is in (not necessarily minimal) bridge
position, and $\alpha$ is contained in a level sphere $S$ as in Case~(1) of
Figure~12 of~\cite{GST}; that is, in the trivial $n$-strand tangle in the
ball in $S^3$ lying above $S$, $\alpha$ connects the endpoints of two
different strands, and the $n$ strands are parallel to a collection of $n$
disjoint arcs in $S$ that are disjoint from the interior of $\alpha$. We
call this a \textit{level arc position} of $\tau$, keeping in mind that
$K_\tau$ is not assumed to be in minimal bridge position. Note, however,
that according to lemma~\ref{lem:GSTlemma}, results of Goda, Scharlemann,
and Thompson show that for a regular tunnel, we may always choose a level
arc position for $\tau$ in which the number of bridges of $K_\tau$
equals~$\br(K_\tau)$.
\end{definition}

\begin{figure}
\begin{center}
\includegraphics[width=\textwidth]{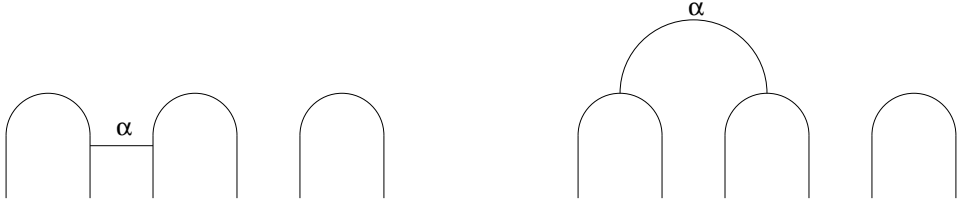}
\caption{Relative bridge counts}
\label{fig:arches}
\end{center}
\end{figure}
\begin{definition}
Fix a level arc position for $\tau$. We refer to the number of bridges of
$K_\tau$ as the \textit{bridge count} of $K_\tau$ for this position, and
denote it by $\bc(K_\tau)$. If $\alpha$ were moved from level arc position
to the position shown in the right-hand diagram of figure~\ref{fig:arches},
each of $K_\lambda$ and $K_\rho$ would contain a certain number of local
maximum points, with the local maximum that lies on $\alpha$ shared by
both. We call those numbers the \textit{relative bridge counts} of
$K_\lambda$ and $K_\rho$ for the level arc position of $\alpha$, and denote
them by $\rbc(K_\lambda)$ and $\rbc(K_\rho)$. Notice that
$\bc(K_\tau)=\rbc(K_\lambda)+\rbc(K_\rho)$.
\end{definition}

\begin{proposition}
Suppose that $\tau$ is in a level arc position, and that a cabling
operation as in figure~\ref{fig:minraise} is performed, producing a new
tunnel $\tau'$ with principal vertex $\{\rho, \tau, \tau'\}$, and producing
a tunnel arc $\sigma$ for which $\tau'$ is in a level arc position. In
particular, $K_{\tau'}\cup \sigma$ contains knots $K_\rho'$ and $K_\tau'$
equivalent to $K_\rho$ and $K_\tau$. Then
\begin{enumerate}
\item $\rbc(K_\rho')=\rbc(K_\rho)$.
\item $\rbc(K_\tau')=\bc(K_\tau)$.
\item $\bc(K_{\tau'}) = \bc(K_\tau) + \rbc(K_\rho)$.
\end{enumerate}
\label{prop:bridge_counts}
\end{proposition}
\begin{proof}
Careful examination of figure~\ref{fig:minraise} verifies the proposition
in case the arc of the original $K_\tau\cup \alpha$ dual to $\lambda$ has
one end that leaves $\alpha$ in the upward direction and one end that
leaves it in the downward direction. There are two other possibilities,
either both ends leave in the upward direction, or both leave in the
downward direction. Very similar constructions verify the proposition in
those two cases.
\end{proof}

The next two results follow easily from
proposition~\ref{prop:bridge_counts}.
\begin{theorem} Let $\tau$ be a regular tunnel. In the principal path 
of $\tau$, in which $\tau=\tau_n$, let $\tau_m$ be the first tunnel of
depth $2$, with principal vertex $\{\tau_{m-2},\tau_{m-1},\tau_m\}$. Put
$\tau_m$ in a level
arc position, and let $u_{m-2}=\rbc(K_{\tau_{m-2}})$ and
$u_{m-1}=\rbc(K_{\tau_{m-1}})$. To the vertices $\tau_m,\ldots\,$,
$\tau_n$, assign values $u_k$ inductively by the rule $u_k=u_i+u_j$, where
the principal vertex of $\tau_k$ is $\{\tau_i, \tau_j,\tau_k\}$. Then the
bridge number of $K_\tau$ is at most~$u_n$.
\label{thm:upper_bound}
\end{theorem}

\begin{corollary} Let $\tau$ be a regular tunnel. In the principal path 
of $\tau$, in which $\tau=\tau_n$, let $\tau_m$ be the first tunnel of
depth $2$, with principal vertex $\{\tau_{m-2},\tau_{m-1},\tau_m\}$.
Suppose that the tunnel $\tau_m$ can be put in a level arc position so that
$\rbc(K_{\tau_j})=\br(K_{\tau_j})$ for $m-2\leq j\leq m-1$, and
$\bc(K_{\tau_m})=\br(K_{\tau_m})$. Then the bridge number of $K_\tau$
equals the value~$c_n$ of theorem~\ref{thm:bridge_numbers}.\par
\label{coro:exact_bridge_number}
\end{corollary}

\begin{figure}
\begin{center}
\includegraphics[width=\textwidth]{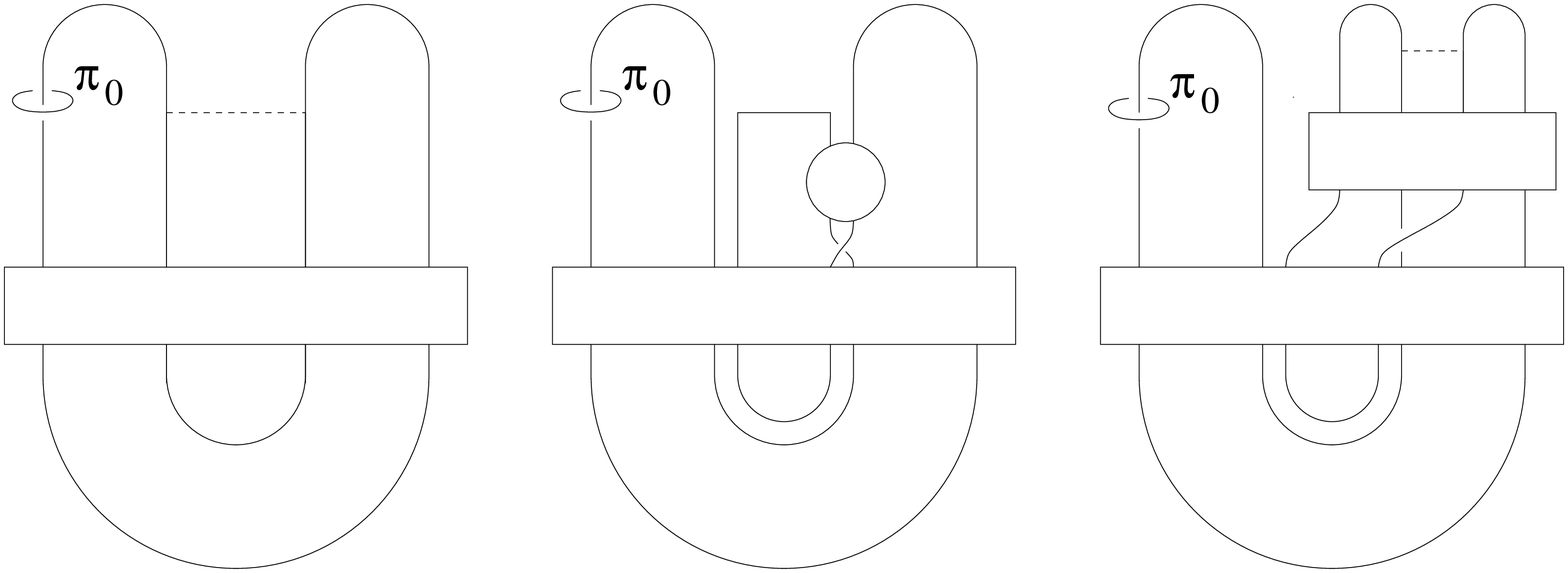}
\caption{Conservative cabling of semisimple tunnels}
\label{fig:semisimple}
\end{center}
\end{figure}
We now focus on cablings that produce semisimple tunnels. When
Proposition~\ref{prop:bridge_counts} is applied at each step of the cabling
sequence of a semisimple tunnel, each cabling construction can be performed
so that the bridge number of the resulting $K_\tau$ is $1$ larger than the
bridge number of the previous knot. Figure~\ref{fig:semisimple} (which, as
in the proof of proposition~\ref{prop:bridge_counts}, admits two variants)
illustrates the inductive process. Each rectangle in that figure represents
a pure $4$-strand braid, and the circle represents a rational tangle. The
first diagram shows a level simple tunnel of a $2$-bridge knot
$K_{\tau_0}$. The next cabling operation is performed, producing a knot
$K_{\tau_1}$ with a semisimple tunnel, as in the second diagram. This is
moved by isotopy to the position in the third diagram; the tunnel is in a
level arc position, and $K_{\pi_0}$ has only one (relative) bridge. The
$k^{th}$ repetition of this construction produces a level arc position of
$\tau_k$ for which $\bc(K)_{\tau_k}=k+2$, $\rbc(K_{\tau_{k-1}})=k+1$, and
$\rbc(K_{\pi_0})=1$. Therefore we have:
\begin{theorem}
Let $\tau$ be a semisimple tunnel produced by $m$ cabling constructions,
and let $\{\pi_0,\rho,\tau\}$ be its principal vertex. Then $\tau$ can be
placed in level arc position so that $\bc(K_\tau)=m+1$, $\rbc(K_\rho)=m$,
and $\rbc(K_{\pi_0})=1$.\par
\label{thm:level_semisimple}
\end{theorem}

\begin{corollary} Suppose that a semisimple tunnel $\tau$ is produced by 
$m$ cabling constructions. The $\br(K_\tau)\leq m+1$.
\end{corollary}

Combining theorem~\ref{thm:level_semisimple} with
theorem~\ref{thm:upper_bound}, we have an upper bound for bridge number:
\begin{theorem} Let $\tau$ be a regular tunnel. In the principal path 
of $\tau$, in which $\tau=\tau_n$, let $\tau_m$ be the first tunnel of
depth $2$. Put $u_{m-2}=m$ and $u_{m-1}=m+1$. To the vertices
$\tau_m,\ldots\,$, $\tau_n$, assign values $u_k$ inductively by the rule
$u_k=u_i+u_j$, where the principal vertex of $\tau_k$ is $\{\tau_i,
\tau_j,\tau_k\}$. Then the bridge number of $K_{\tau_n}$ is at most~$u_n$.
\label{thm:general_upper_bound}
\end{theorem}

We can give a universal upper bound for the bridge number of $K_\tau$ in
terms of the number of cablings that produce~$\tau$.
\begin{figure}
\begin{center}
\includegraphics[height=30 ex]{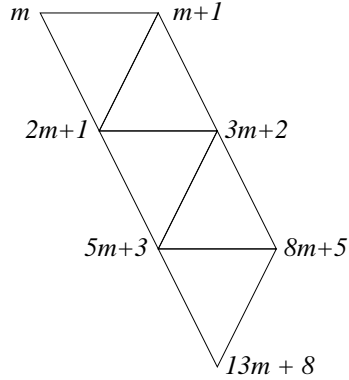}
\caption{The fastest growing upper bounds for bridge number, starting with
the last two semisimple tunnels in the cabling sequence.}
\label{fig:universal}
\end{center}
\end{figure}
We denote by $(F_1,F_2,\ldots)$ the Fibonacci sequence
$(1,1,2,3,\ldots)$.
\begin{theorem} 
Let $\tau$ be a regular tunnel produced by $n$ cabling operations, of
which the first $m$ produce semisimple tunnels. Then
$\br(K_{\tau_n})\leq mF_{n-m+2}+F_{n-m+1}$.
\label{thm:ub}
\end{theorem}
\begin{proof}
Figure~\ref{fig:universal} shows the type of principal path for
which the $u_k$ in theorem~\ref{thm:general_upper_bound} grow most rapidly
(in figure~\ref{fig:universal}, the top two vertices are $\tau_{m-2}$ and
$\tau_{m-1}$, the last two semisimple tunnels that appear in the cabling
sequence of $\tau_n$). Putting $u_{m-2}=m$, $u_{m-1}=m+1$, and
$u_k=u_{k-1}+u_{k-2}$ for $k\geq m$, theorem~\ref{thm:general_upper_bound}
gives $\br(K_{\tau_n})\leq u_n$. Since $u_{m-2}=m\cdot F_1$ and
$u_{m-1}=m\cdot F_2 + F_1$, the $n-m$ additional recursions give the
estimate in the theorem.
\end{proof}

For a fixed $n$, the largest upper bound in theorem~\ref{thm:ub} occurs
when $m=2$. We finish by showing that theorem~\ref{thm:ub} is sharp for
that case.
\begin{theorem} The maximum bridge number of any tunnel number
$1$-knot having a tunnel produced by $n$ cabling operations
is~$F_{n+2}$.\par
\label{thm:max}
\end{theorem}
\begin{proof}
Since the minimum possible value for $m$ in theorem~\ref{thm:ub} is $2$,
any tunnel $\tau$ produced by $n$ cabling operations has $\br(K_\tau)\leq
2F_n+F_{n-1}=F_{n+2}$. Now, let $\tau_0$ be any simple tunnel.
In~\cite{CM}, the slope sequences for the semisimple tunnels of $2$-bridge
knots were calculated, in particular finding that each cabling had slope of
the form $\pm 2 + 1/k$ for some integer $k$. Perform any cabling on
$K_{\tau_0}$ whose slope is not of this form, to produce a tunnel $\tau_1$
for which $K_{\tau_1}=3$. We now perform cabling constructions following
the principal path indicated in figure~\ref{fig:universal} with
$m=2$. Theorem~\ref{thm:bridge_numbers} shows that after $n$ cabling
constructions, $\br(K_\tau)\geq F_{n+2}$.
\end{proof}

\section{Tunnels of torus knots}
\label{sec:torus_knots}

The tunnels of torus knots were analyzed by M. Boileau, M. Rost, and
H. Zieschang \cite{B-R-Z}. There are two $(1,1)$-tunnels, and a third
``short'' tunnel represented by an arc that cuts straight across the
complementary annulus when the knot is regarded as being contained in a
standard torus. In certain cases, some of these tunnels are equivalent. In
this section, we will analyze the cabling sequences for the short
tunnels. In particular, we will see that their depths are arbitrarily
large. On the other hand, all have distance~$2$, as we will now
verify, while setting some notation for this section.

Consider a (nontrivial) $(p,q)$ torus knot, contained in a standard torus
$T$ in $S^3$, bounding a solid torus $W\subset \R^3\subset S^3$. The short
tunnel is represented by an arc in $T$. Let $H$ be a regular neighborhood
of the knot and tunnel arc, chosen to intersect $T$ in a regular
neigbhorhood of the knot and tunnel arc in $T$. Let $\tau$ be the cocore
disk of the short tunnel, so that the torus knot is $K_\tau$. Now $K_\tau$
is isotopic to a loop $C$ in $\partial H$ that lies entirely outside of $W$
and is disjoint from $\partial \tau$. But $C$ is also disjoint from the
disk $\overline{S^3-H}\cap T$, showing that $\tau$ has distance~$2$.

\begin{figure}
\begin{center}
\includegraphics[width=\textwidth]{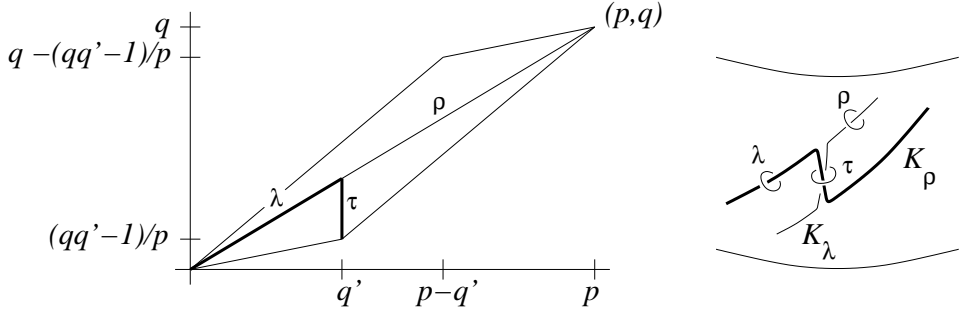}
\caption{The properties of $q'$. The darker segments correspond to
$K_\rho$, a $(q',(qq'-1)/p)$ torus knot. The picture on the right shows
$K_\lambda$ in the torus $T\subset S^3$, and $K_\rho$ pulled slightly
outside of $T$.}
\label{fig:qprime}
\end{center}
\end{figure}
We now begin our calculation of the slope invariants of the short tunnels.
First we examine a cabling operation that takes a short tunnel $\tau$ and
produces a short tunnel of a new torus knot.

Assume for now that both $p$ and $q$ are positive. Since the $(p,q)$ and
$(q,p)$ torus knots are isotopic, we may further assume that~$p>q$. Let
$q'$ be the integer with $0<q'<p$ such that $qq'\equiv
1\pmod{p}$. Figure~\ref{fig:qprime} illustrates the features of $q'$. If
the principal pair $\{\lambda,\rho\}$ of $\tau$ is positioned as shown in
figure~\ref{fig:qprime} (our inductive construction of these tunnels will
show that the pair shown in the figure is indeed the principal pair), then
$K_\rho$ is a $(q',(qq'-1)/p)$ torus knot, and $K_\lambda$ is a $(p-q',
q-(qq'-1)/p)$ torus knot. We set $(p_1,q_1)= (q',(qq'-1)/p)$ and $(p_2,q_2)
= (p-q', q-(qq'-1)/p)$, so that $K_\rho$ and $K_\lambda$ are respectively
the $(p_1,q_1)$ and $(p_2,q_2)$ torus knots.

In figure~\ref{fig:qprime}, the linking number of $K_\rho$ with
$K_\lambda$, up to sign conventions, is $q_1p_2$. One way to see this is to
note that a Seifert surface for $K_\lambda$ can be constructed from $q_2$
meridian disks for $W$ and $p_2$ meridian disks for the ``outside'' solid
torus $\overline{S^3-W}$. When $K_\rho$ is pulled slightly outside of $W$,
as indicated in figure~\ref{fig:qprime}, each of the $p_2$ meridian disks
for the outside solid torus has $q_1$ intersections with
$K_\rho$, all crossing the disks in the same direction.

\begin{figure}
\begin{center}
\includegraphics[width=\textwidth]{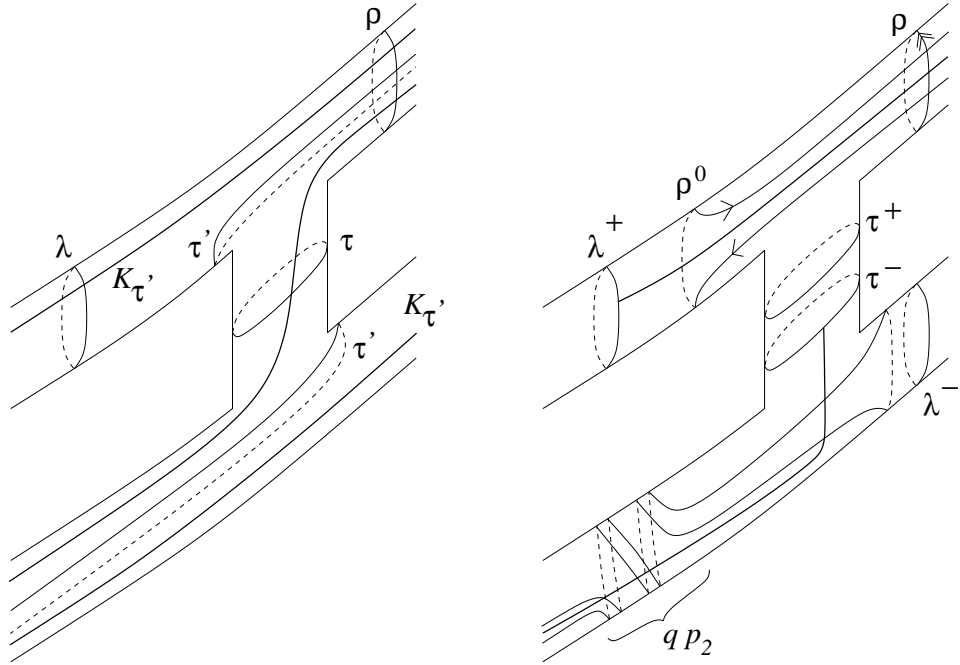}
\caption{The cabling construction that replaces $\rho$
(compare with figure~\ref{fig:slope_coords}). The left drawing
shows the new tunnel disk $\tau'$ and the knot $K_{\tau'}$. The right
drawing shows a cabling arc for $\tau'$, running from $\lambda^+$ to
$\tau^-$, and the disks $\rho$ and $\rho^0$ used to calculate its
slope. The $qp_2$ turns of $\rho^0$, with the case $qp_2=2$ drawn in the
figure, make the copies of $K_\tau$ and $K_\lambda$ in its complement have
linking number $0$.}
\label{fig:torus_tunnel}
\end{center}
\end{figure}
Figure~\ref{fig:torus_tunnel} shows the new tunnel disk $\tau'$ for a
cabling construction that produces a $(p+p_2,q+q_2)$-torus knot
$K_{\tau'}$. This disk meets $T$ perpendicularly. The drawing on the right
in figure~\ref{fig:torus_tunnel} illustrates the setup for the calculation
of the $(\{\lambda,\tau\};\rho)$-slope pair of $\tau'$.
Examination of that drawing shows that the
slope pair of $\tau'$ is~$[1,2qp_2+1]$.

As usual, let $U=\begin{pmatrix}1&1\\0&1\end{pmatrix}$ and
$L=\begin{pmatrix}1&0\\1&1\end{pmatrix}$. If $K_1$ is a $(p_1,q_1)$-torus
knot and $K_2$ is a $(p_2,q_2)$-torus knot, we denote by $M(K_1,K_2)$ the
matrix $\begin{pmatrix} p_1 & q_1 \\ p_2 & q_2\end{pmatrix}$. In our case,
this is the matrix $M(K_\rho, K_\lambda)$. Adding the rows of
$M(K_\rho,K_\lambda)$ gives $(p,q)$, corresponding to $K_\tau$, so
\[ M(K_\tau,K_\lambda) = U\cdot M(K_\rho,K_\lambda)\ .\]
The left drawing of figure~\ref{fig:torus_tunnel} can be repositioned by
isotopy so that $\lambda$, $\tau$, and $\tau'$ look respectively as did
$\lambda$, $\rho$, and $\tau$ in the original picture, with $\tau'$ as the
tunnel of the $(p+p_2,q+q_2)$ torus knot. Thus the procedure can be
repeated, each time multiplying the matrix by another factor of~$U$.

\begin{figure}
\begin{center}
\includegraphics[width=0.45 \textwidth]{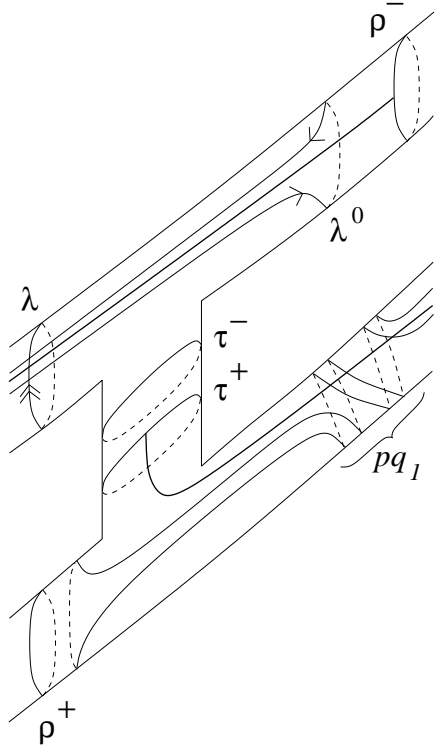}
\caption{Calculation of the slope of $\tau'$ for a cabling construction
replacing $\lambda$. The cabling arc runs from $\rho^-$ to $\tau^+$.}
\label{fig:torus_tunnel1}
\end{center}
\end{figure}
Figure~\ref{fig:torus_tunnel1} shows the calculation of the slope of the
cabling construction replacing $\lambda$ by a new tunnel $\tau'$, with the
effect that 
\[ M(K_\rho, K_\tau) = L\cdot M(K_\rho,K_\lambda)\ .\]
Its slope pair works out to be $[1,2qp_1-1]$. One might expect $2qp_1+1$ as
the second term, in analogy with the construction replacing
$\rho$. However, as seen in figure~\ref{fig:torus_tunnel1}, the $pq_1$
twists needed in $\lambda^0$ are in the same direction as the twists in the
calculation for $\rho$, not in the mirror-image sense. This results in two
fewer crossings of the cabling arc for $\tau'$ with $\lambda^0$ than
before. In fact, the slope pairs for the two constructions can be described
in a uniform way: For either of the matrices $M(K_\tau,K_\lambda)$ and
$M(K_\rho,K_\tau)$, a little bit of arithmetic shows that the second entry
of the slope pair for the cabling operation that produced them is the sum
of the product of the diagonal entries and the product of the off-diagonal
entries, that is, $[1,pq_2+qp_2]$ in the first case and $[1,pq_1+qp_1]$ in
the second.

We can now describe the complete cabling sequence. Still assuming that $p$
and $q$ are both positive and $p>q$, write $p/q$ as $[n_1,n_2,\ldots, n_k]$
with all $n_i$ positive. We may assume that $n_k\neq 1$. According as $k$
is even or odd, consider the product $U^{n_k}L^{n_{k-1}}\cdots
U^{n_2}L^{n_1}$ or $L^{n_k}U^{n_{k-1}}\cdots U^{n_2}L^{n_1}$. Start with a
trivial knot regarded as a $(1,1)$-torus knot, and a tunnel positioned so
that $K_\rho$ is a $(1,0)$-torus knot and $K_\lambda$ is a $(0,1)$-torus
knot. The corresponding matrix $M(K_\rho,K_\lambda)$ is the identity
matrix.  Multiplying by $L^{n_1}$ has the effect of doing $n_1$ trivial
cabling constructions, each with slope $1$, and ending with the trivial
knot positioned as an $(n_1+1,1)$-torus knot.  At that stage,
$M(K_\rho,K_\lambda)$ is $\begin{pmatrix} 1 & 0 \\ n_1 &
1\end{pmatrix}$. Then, multiplying by $U$ corresponds to a true cabling
construction. In the above notation, the new matrix $M(K_\tau,K_\lambda)$
is $\begin{pmatrix} n_1+1 & 1 \\ n_1 & 1\end{pmatrix}$, and the knot
$K_{\tau'}$ is a $(2n_1+1,2)$-torus knot. As explained above, the
construction has slope pair $[1,2n_1+1]$, so the simple slope is
$m_0=[1/(2n_1+1)]$. Continue by multiplying $n_2-1$ additional times by
$U$, then $n_3$ times by $L$ and so on, performing additional cabling
constructions with slopes calculated as above from the matrices of the
current $K_\rho$, $K_\lambda$, and $K_\tau$.

At the end, there is no cabling construction corresponding to the last
factor $L$ or $U$. For specificity, suppose $k$ was even and the product
was $U^{n_k}L^{n_{k-1}}\cdots U^{n_1}L^{n_1}$.  At the last stage, we apply
$n_k-1$ cabling constructions corresponding to multiplications by $U$, and
arrive at a tunnel $\tau$ for which $M(K_\rho,K_\lambda)$ is
$U^{n_k-1}L^{n_{k-1}}\cdots U^{n_2}L^{n_1}$. The sum of the rows is then
$(p,q)$ (multiplying by $U$ and using the case ``$q/s$'' of lemma~14.3
of~\cite{CM}), so $K_\tau$ is the $(p,q)$ torus knot. The case when $k$ is
odd is similar (using the ``$p/r$'' case of lemma~14.3 of~\cite{CM} at the
end). In summary, there are $-1+\sum_{i=2}^{k} n_i$ (nontrivial) cabling
constructions, whose slopes can be calculated as above.

Suppose now that $p$ is positive but $q$ is negative. We may assume that
$p>\vert\,q\,\vert$. We have already found the cabling sequence for the
case of the $(p,-q)$-torus knot, and since reversing orientation negates
the slope parameters in the Parameterization Theorem (remark~12.5
of~\cite{CM}), we need only negate its slopes to obtain the cabling
sequence for the $(p,q)$-torus knot.

\begin{figure}
\begin{center}
\includegraphics[height=55 ex]{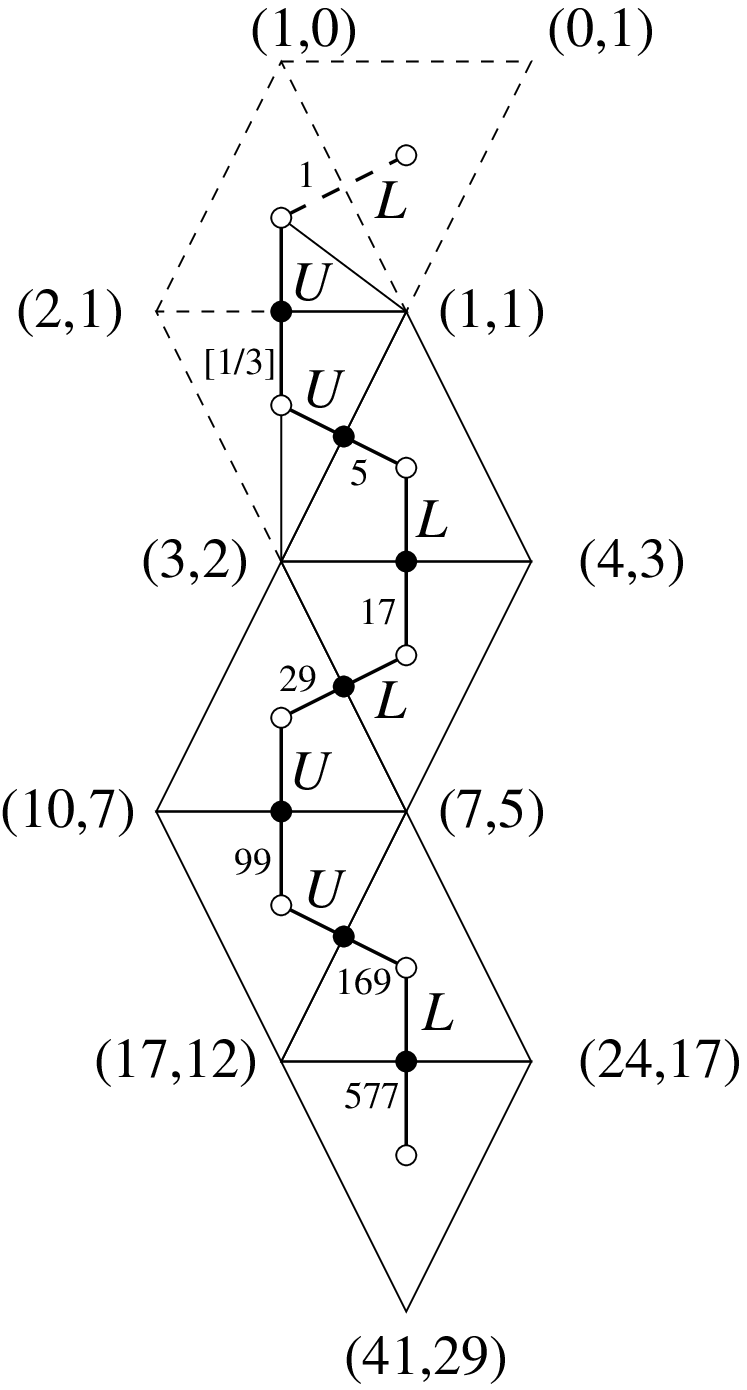}
\caption{Slowest growth of bridge number as a function of depth for torus
knot tunnels, corresponding to the continued fraction expansion
$41/29=[1,2,2,2,2]$. The $(41,29)$ torus knot has the smallest bridge
number of any torus knot with a depth~$4$ tunnel.}
\label{fig:cheapest}
\end{center}
\end{figure}
\begin{remark}\label{rem:torus_tunnel_bridge_number}
Figure~\ref{fig:cheapest} shows an initial segment of the principal path
for the tunnels of the $(p,q)$-torus knots for the continued fractions
$p/q=[1,2,2,\ldots,2]$, which limit to $\sqrt{2}$. Notice that this is the
path of cheapest descent from figure~\ref{fig:descent}. The small numbers
along the path are the slopes, the letters indicate whether the
constructions correspond to multiplication by $U$ or by $L$, and the pairs
show the $(p,q)$ for the torus knots determined by the tunnel at each step.
The first nontrivial cabling, with $m_0=[1/3]$, produces a $(3,2)$-torus
knot, and the second produces a $(4,3)$-torus knot with bridge number
$3$. Since we always have $p>q$, the bridge number is simply the value of
$q$. These obey the recursion of corollary~\ref{coro:bridge_numbers},
starting with $b_2=2$ and $b_3=3$. Since the cabling sequence for the short
tunnel $\tau$ of any torus knot contains only one two-bridge knot (the
$(2n_1+1,2)$-torus knot produced by the first nontrivial cabling), there is
no regular torus knot tunnel which has $b_3=2$. Therefore each $b_{2d}$ in
this sequence (which also occur for its reversed-orientation sequence,
where the negatives of these slopes are used) gives the minimum bridge
number for a torus knot with a tunnel of depth~$d$. Finally, we note that
$b_{2d}=a_d$ where $a_1=2$, $a_2=5$, and $a_d=2a_{d-1}+a_{d-2}$ for $d\geq
3$. The asymptotic growth rate of this sequence is proportional to
$(1+\sqrt{2})^d$, which as we have seen is the minimum rate in general.
\end{remark}

\begin{remark}
Considering the preceding example, we can see how to determine the depth
from the continued fraction expansion $[n_1,n_2,\ldots,n_k]$ of
$p/q$. Ignore $n_1$, since it corresponds to cablings which produce
the trivial tunnel. Basically, each of the remaining $n_i$ increases the
depth by $1$. However, for a block $[n_i,\ldots, n_{i+j}]$ with each
$n_\ell=1$ for $i\leq \ell \leq i+j$, and $n_{i-1}\neq 1$ and
$n_{i+j+1}\neq 1$, only the cablings corresponding $n_i$, $n_{i+2}$,
$n_{i+4}$, and so on increase the depth (consider the principal path drawn
as in figure~\ref{fig:corridor}). Also, a final block $[n_{k-1},n_k]=[1,2]$
increases the depth by only $1$, since the final $2$ represents only a
single cabling operation.
\end{remark}

The exceptional cases of M. Boileau, M. Rost, and H. Zieschang
\cite{B-R-Z}, that is, the cases when there are fewer than three tunnels,
are exactly the cases when $\tau$ is semisimple. To understand this
computationally, we may assume that $p>q\geq 2$, and we find:
\smallskip

\noindent \textsl{Case I:} $p\equiv 1\pmod q$.  We have $p/q=[n_1,q]$ and
we examine $U^qL^{n_1}$. There are $n_1$ trivial cablings, producing the
trivial tunnel, then there are $q-1$ cabling constructions retaining one of
the original arcs of the trivial knot, showing that $\tau$ is
semisimple.\par
\smallskip

\noindent \textsl{Case II:} $p\equiv -1\pmod q$.

In these cases, $p/q$ is $[n_1,1,q]$, with $q>1$. Examining
$L^{q-1}UL^n$, the first nontrivial cabling corresponds to $U$, and
produces a simple tunnel, then the $q-1$ cablings corresponding to
the $L^q$ retain one of the original arcs of $K_{\tau_0}$. Thus these
are also semisimple tunnels.
\smallskip

In all other cases, the continued fraction expansion of $p/q$ either has
more than three terms, or has second term greater than $1$, so the regular
tunnel is regular. We have verified the equivalence of the first two
conditions in the following proposition. As already noted, the equivalence
of the first and third is from~\cite{B-R-Z}.

\begin{proposition} For the
$(p,q)$ torus knot $K_{(p,q)}$, every tunnel has distance~$2$. The
following are equivalent:
\begin{enumerate}
\item $p\not\equiv \pm1 \pmod q$.
\item The short tunnel is regular.
\item $K_{(p,q)}$ has exactly three tunnels.
\end{enumerate}
\label{prop:torus_regular}
\end{proposition}

We have implemented the algorithms for the slope sequence and the depth of
the short tunnel computationally~\cite{slopes}. Some sample calculations
are:
\smallskip

\noindent \texttt{TorusKnots$>$ slopes 41 29}\\
\texttt{[ 1/3 ], 5, 17, 29, 99, 169, 577}
\smallskip

\noindent \texttt{TorusKnots$>$ slopes 181 (-48)}\\
\texttt{[ 6/7 ], -15, -23, -31, -151, -271, -883, -2157, -3431}
\smallskip

\noindent \texttt{TorusKnots$>$ depth 41 29}\\
\texttt{4}
\smallskip

\noindent \texttt{TorusKnots$>$ [ depth 41 n | n <- [2..40] ]}\\
\texttt{[1,1,1,1,1,1,1,2,1,2,3,2,1,2,3,3,3,2,1,1,2,3,3,3,3,2,3,4,3,2,3,
2,2,2,2,2,2,2,1]}
\smallskip

\noindent The last command produces a list of the depths of the short
tunnels for the torus knots $K_{(41,2)}$ through $K_{(41,40)}$.

\section{The case of tunnel number $1$ links}
\label{sec:links}

As explained in ~\cite{CM}, our entire theory can be adapted to include
tunnels of tunnel number $1$ links simply by adding the separating disks as
possible slope disks. The full disk complex $\K(H)$ is only slightly more
complicated than $\D(H)$. Each separating disk is disjoint from only two
other disks, both nonseparating, so the additional vertices appear in
$2$-simplices attached to $\D(H)$ along the edge opposite the vertex that
is a separating disk. The quotient $\K(H)/\G$ only has such additional
$2$-simplices: (1)~there is a unique orbit of ``primitive'' separating
disk, consisting of separating disks disjoint from a primitive pair, which
are exactly the intersections of splitting spheres with $H$. Their orbit
$\sigma_0$ is a vertex of a ``half-simplex'' $\langle \sigma_0, \pi_0,
\mu_0\rangle$ attached to $\D(H)/\G$ along $\langle \pi_0,
\mu_0\rangle$. It is the unique tunnel of the trivial $2$-component link,
and has simple slope $\infty$. (2)~Simple separating disks lie in
half-simplices attached along $\langle \pi_0, \mu_0\rangle$, just like
nonseparating simple disks. Their simple slopes are $[p/q]$ with $q$
even. (3)~the remaining separating disks lie in $2$-simplices attached
along edges of $\D(H)/\G$ spanned by two (orbits of) disks, at least one of
which is nonprimitive. A single ``Y'' is added to $\T$ for each added
$2$-simplex (or a folded ``Y'', for the half-simplices). The link in
$\K'(H)/\G$ of a link tunnel is simply the top edges (or top edge, for
the trivial and simple tunnels) of such a~``Y''.

Cabling operations differ only in allowing separating slope disks, which
produce a tunnel of a tunnel number~$1$ link.  The cabling sequence ends
with the first separating slope disk, and cannot be continued.  The
Parameterization Theorem holds as stated, except allowing $q_n$ to be even,
and allowing $m_0=\infty$ for the unique tunnel of the trivial link.

For link tunnels, the distance and depth invariants are defined as for knot
tunnels. Depth~$1$ tunnels are the tunnels of links with one component
unknotted. The other component must be a $(1,1)$-knot, and the link must
have torus bridge number~$2$~\cite{CM}.  Lemma~\ref{lem:JohnsonDistance}
holds when $\tau$ is separating, in fact the argument is an easier version
of the argument in~\cite{JohnsonBridgeNumber}, so
theorem~\ref{thm:uniquetunnel} and corollary~\ref{coro:amphichiral} hold
for links as well as knots.

For a tunnel $\sigma$ of a tunnel number $1$ knot, there is a version of a
giant step that produces a tunnel number~$1$ link.  Choose any loop
$L_1$ in $\partial H$ that crosses $\sigma$ exactly once, and let $\sigma'$
be the frontier of a regular neighborhood of $L_1\cup \sigma$ in $H$. Since
$\sigma'$ is separating, the core circles of its complementary solid tori
form a tunnel number~$1$ link with tunnel $\sigma'$; one of these core
circles is isotopic to $L_1$.  One might even describe a giant step
starting from a tunnel of a link, but this is of little interest since such
a giant step could not appear in a minimal giant step sequence starting
from~$\pi_0$, because the two disks disjoint from a separating disk are
also disjoint from each other. Section~\ref{sec:GST_move_seqs} adapts
almost word-for-word to allow tunnels of links.

The proof of theorem~\ref{thm:bridge_numbers} adapts without difficulty to
the case of links since the tunnel leveling of~\cite{GST} applies to links
as well as knots. Since the geometric constructions in
sections~\ref{sec:bridge_number_growth} and~\ref{sec:bridgenum_conj} also
work for cabling constructions that produce links, the results of both
those sections apply just as well to links.

\bibliographystyle{amsplain}

\begin{thebibliography}{10}
\bibitem{Akbas} E. Akbas, A presentation of the automorphisms of the
$3$-sphere that preserve a genus two Heegaard splitting, Mathematics ArXiv
math.GT/0504519.
\end{thebibliography}

\end{document}